\documentclass{article}
\usepackage[utf8]{inputenc}
\usepackage{graphicx}
\usepackage{tikz-cd}
\usepackage{grffile}
\usepackage{ctable}
\usepackage{appendix}
\usepackage{geometry}
\usepackage{float}
\usepackage[footnotesize]{caption}
\usepackage{subcaption}
\usepackage{epstopdf}

\usepackage{amssymb, nccmath, comment}
\usetikzlibrary{matrix,calc}
\usepackage{amsmath}
\usepackage{mathrsfs}
\usepackage[dvips]{epsfig}

\usepackage[english]{babel} 
\usepackage{enumerate}
\usepackage[utf8]{inputenc}
\usepackage[all]{xy}

\usepackage [autostyle, english = american]{csquotes}
\usepackage{bm}
\MakeOuterQuote{"}
\newcommand{\Z}{\mathbb{Z}}
\newcommand{\N}{\mathbb{N}}

\newcommand{\C}{\mathbb{C}}
\newcommand{\R}{\mathbb{R}}

\newcommand{\half}{\frac{1}{2}}

\newcommand{\ep}{\varepsilon}

\newcommand{\T}{\mathcal{T}}

\newcommand{\wt}{\widetilde}

\newcommand{\ran}{\mathrm{Ran\,}}

\DeclareMathOperator{\sech}{sech}
\newcommand*{\bigchi}{\mbox{\Large$\chi$}}
\newcommand{\tbeta}{\Tilde{\beta}}
\newcommand{\tc}{\Tilde{c}}

\definecolor{Green}{rgb}{0.,0.4,0.}

\renewcommand{\Re}{\mathrm{Re}}
\renewcommand{\Im}{\mathrm{Im}}

\newcommand{\Rmnum}[1]{\uppercase\expandafter{\romannumeral #1\relax}}

\DeclareFontFamily{OT1}{pzc}{}

\DeclareMathAlphabet{\mathpzc}{OT1}{pzc}{m}{it}

\newtheorem{Lemma}{Lemma}[section]
\newtheorem{Theorem}{Theorem}
\newtheorem{Proposition}[Lemma]{Proposition}

\newtheorem{Remark}[Lemma]{Remark}

\newenvironment{Proof}[1][\unskip]%
 {\begin{trivlist} \item[]{\bf Proof #1. }}%
 {\hspace*{\fill}$\rule{.4\baselineskip}{.4\baselineskip}$\end{trivlist}}

 {\begin{trivlist}\item[]\textbf{Acknowledgments.}}{\end{trivlist}}

 {\begin{center}\textbf{Abstract}}{\end{center}}

\makeatletter\@addtoreset{figure}{section}\makeatother

\makeatletter \@addtoreset{equation}{section} \makeatother

\makeatletter
\newsavebox{\@brx}
\newcommand{\llangle}[1][]{\savebox{\@brx}{\(\m@th{#1\langle}\)}%
  \mathopen{\copy\@brx\kern-0.5\wd\@brx\usebox{\@brx}}}
\newcommand{\rrangle}[1][]{\savebox{\@brx}{\(\m@th{#1\rangle}\)}%
  \mathclose{\copy\@brx\kern-0.5\wd\@brx\usebox{\@brx}}}
\makeatother

\definecolor{Green}{rgb}{0.,0.4,0.}

\renewcommand{\leq}{\leqslant}
\renewcommand{\geq}{\geqslant}

\def\XXint#1#2#3{{\setbox0=\hbox{$#1{#2#3}{\int}$}
     \vcenter{\hbox{$#2#3$}}\kern-.5\wd0}}

\newfam\bifam
\font\tenbi=cmmib10 scaled \magstep1 \font\sevenbi=cmmib10 at 11pt
\font\fivebi=cmmib10 at 6pt \textfont\bifam = \tenbi
\scriptfont\bifam = \sevenbi \scriptscriptfont\bifam= \fivebi

\title{Shifting consensus in a biased compromise model}
\begin{document}
\begin{center}
{\fontsize{15}{15}\fontfamily{cmr}\fontseries{b}\selectfont{Shifting consensus in a biased compromise model}}\\[0.2in]
Olivia Cannon$^2$, Ty Bondurant$^3$, Malindi Whyte$^{4}$, and Arnd Scheel$^{2,}$\footnote{The authors acknowledge partial support through  grants NSF DMS-1907391 and NSF DMS-2205663.}\\[0.1in]
\textit{\footnotesize 
$^2$ University of Minnesota, School of Mathematics,   206 Church St. S.E., Minneapolis, MN 55455, USA} \\
\textit{\footnotesize 
$^3$ Mathematics Department, Georgia Institute of Technology, North Ave NW, Atlanta, GA 30332 USA} \\ 
\textit{\footnotesize $^4$
Department of Mathematics, Wake Forest University, 127 Manchester Hall,
Winston-Salem, NC 27109, USA}
\end{center}

\begin{abstract}
    We investigate the effect of bias on the formation and dynamics of political parties in the bounded confidence model. For weak bias, we quantify the change in average opinion and potential dispersion and decrease in party size. For nonlinear bias modeling self-incitement, we establish coherent drifting motion of parties on a background of uniform opinion distribution for biases below a critical threshold where parties dissolve. Technically, we use geometric singular perturbation theory to derive drift speeds, we rely on a nonlocal center manifold analysis to construct drifting parties near threshold, and we implement numerical continuation in a forward-backward delay equation to connect asymptotic regimes.  
\end{abstract}

\section{Introduction}

The bounded confidence model \cite{HK2002} has been pivotal in the study of social dynamics, providing a mechanism for the formation of opinion clusters, also called parties. Agents are attributed numerical values of opinions. They interact and change their opinions through compromise, but only with other agents whose opinions are sufficiently close, i.e. within a bounded confidence interval. The model can be framed in many different formulations -- stochastic Markov processes, deterministic mean-field equations, discrete or continuous opinion values -- but the qualitative phenomena are similar: a uniform distribution of opinions is an unstable steady-state, and fluctuations lead to the formation of clusters, often regularly spaced \cite{Ben_Naim_2005,bennaimscheel2016, blondel2009statedependent,fortunato2005vector,HK2002, lorenz2007survey, lorenz2007stabilization}. 

This model has been used in large part to study mechanisms of polarization, and to that end, many modifications have been made, including for instance variations of the confidence interval between agents, introduction of a small number of agents who do not compromise ('stubborn' agents), and variations in the probability of interaction \cite{carletti2007propaganda, chen2016opinion, delvicario2016, douven2020misdisinformation, HAN2019121791, kertesz2019algorithmic,  mathias2016, Quattrociocchi2014OpinionDO, shensun2009}. 
The present work is concerned with  \emph{drift} of opinion clusters; that is, with the continuous movement of clusters toward one extreme of the opinion spectrum. While drift caused by asymmetric confidence has been reported in \cite[\S4.2]{HK2002}, there appear to be few systematic computational or analytical studies of this phenomenon. We are interested in this effect as a self-organized phenomenon, caused by behavioral bias in individual agents. We model this through the addition of bias terms to the bounded confidence model. We will demonstrate how bias terms typically lead to drift of average opinions in a party, but may also lead to disintegration of individual  parties (see \S{} \ref{s:nonincite}). Our focus therefore will be on a nonlinear quadratic bias term that corresponds to self-incitement and which leads to both persistent and coherent  drift, avoiding in particular a dispersal and disintegration of the party. The focus on this quadratic bias is also rooted in its relevance as a model for the common sociological phenomenon of group polarization, in which the size of a group is related to the strength of the push to adopt more extreme opinions \cite{group_polarization}. 



To be specific,  we study the following deterministic mean-field model for the evolution of populations of agents $P_n$ with opinion $n\in\Z$, 
\begin{equation} \frac{dP_n}{dt} = 2P_{n+1}P_{n-1} - P_n(P_{n+2} + P_{n-2}) + \beta(P^2_{n+1}-P_n^2), \qquad n \in \Z.\end{equation} 
Here, $\beta=0$ corresponds to the deterministic Hegselmann-Krause bounded-confidence model on a lattice, where populations $P_{n-1}$ and $P_{n+1}$ interact with mass-action rates to form opinions $P_n$. We added the self-incitement bias term $\beta(P^2_{n+1}-P_n^2)$ which can be interpreted as that an individual agent decreases their opinion value by one with probability proportional to the population size with the same opinion: interactions between "same-opinion agents" leads to opinion drift toward the extreme. The strength of the bias term is encoded in the parameter $\beta > 0$. 

Our results for this model can roughly be summarized as follows:
\begin{enumerate}
\item For supercritical, strong bias, $\beta>2$, formation of political parties is suppressed and uniform distribution of opinions is stable; 
\item For weak bias, drift speeds $c$ are at leading order proportional to bias and party mass, $c\sim \frac{2\beta m}{\pi}$;
\item For subcritical bias $\beta\lesssim 2$, we establish rigorously  coherent party drift on a constant background of size $m$ with speed $c\sim 4m$;
\item For bias $0<\beta<2$, we find drifting parties by numerical continuation and find that the background population that supports coherent party drift is exponentially small, $\exp(-const/\beta)$.
\end{enumerate}
The results in (i)-(iii) are analytical. Only in the regime (iii) are we able to establish existence of coherent party drift, while (ii) leaves open the possibility of eventual  dispersal of a party. We also briefly discuss some intriguing phenomena related to stability and instability of drifting parties. 

Technically, the results in (ii) rely on a leading order computation of a flow on a slow manifold using geometric singular perturbation theory, while (iii) uses a recently introduced novel method for analyzing coherent structures in nonlocally coupled equations.

\paragraph{Outline.} In \S{}\ref{s:defs}, we collect information about the bounded confidence model and the model with bias terms added, as well as the linearization at single-party and uniform states. We establish in \S{}\ref{s:small} the speed of drift for $\beta \ll 1$, case (ii), using methods from geometric singular perturbation theory. The $\beta \nearrow 2$ regime, case (iii), is treated in \S{}\ref{s:large}, where drifting solutions are established using nonlocal methods. In \S{}\ref{s:num}, we describe numerical approaches and results, in particular concerning case (iv), above. 
Section \S{}\ref{s:nonincite} describes numerical evidence for lack of coherence when the equation is posed with other bias terms.  
We conclude with a brief discussion.







\section{The bounded confidence model: equilibria, stability, and bias}\label{s:defs}

The bounded confidence model \begin{equation}\label{e:bc} \frac{dP_n}{dt} = 2P_{n+1}P_{n-1} - P_n(P_{n+2} + P_{n-2}), \qquad P_n\geq 0, \ n\in\Z, \end{equation} describes the dynamics of an opinion distribution $P_n(t),n\in\Z$, where $P_n(t)$ represents the local population size with opinion $n$ at time $t$. Agents at sites $n - 1$ and $n+1$ compromise through interaction, moving to opinion $n$, while agents at site $n$ compromise with those at $n +2$ and $n-2$, leaving site $n$. There is no interaction between agents at a distance greater than 2. The equation clearly perserves mass and average opinion (or first moment),
\[
\frac{d}{dt} \sum_n P_n=0, \qquad 
\frac{d}{dt} \sum_n nP_n=0. 
\]
provided reasonable conditions such as sufficient localization at $|n|\to\infty$.
The dynamics of \eqref{e:bc} are  to some extent understood. The consensus process leads to the formation of opinion clusters, which one refers to as parties. For \eqref{e:bc}, these parties can be supported on 1 or 2 opinions sites, that is, for instance, $P_0=1, P_j=0$ for $j\neq 0$, or $P_0=\alpha$, $P_1=1-\alpha$, $P_j=0$ for $j\not\in \{0,1\}$. Parties  separated by at least 2 empty sites between them will not interact and one can in this fashion generate equilibrium states with multiple parties. Political parties, once well established after initial transients, appear to be very robust, although introduction of agents away from existing  parties may lead to the formation of new parties. 
In addition to these single-party states, supported on one or two opinion sites, and the associated well separated multi-party states, the equation also supports equilibria of uniform opinion distribution, $P_n\equiv m$ for all $m$.  We refer to \cite{Ben_Naim_2005} for background on the model and its dynamics. We note however that many questions of stability have not been answered in a precise mathematical fashion, possibly due to the abundance and complexity of possible equilibrium configurations. 

\begin{figure}
    \centering
\includegraphics[scale=0.35]{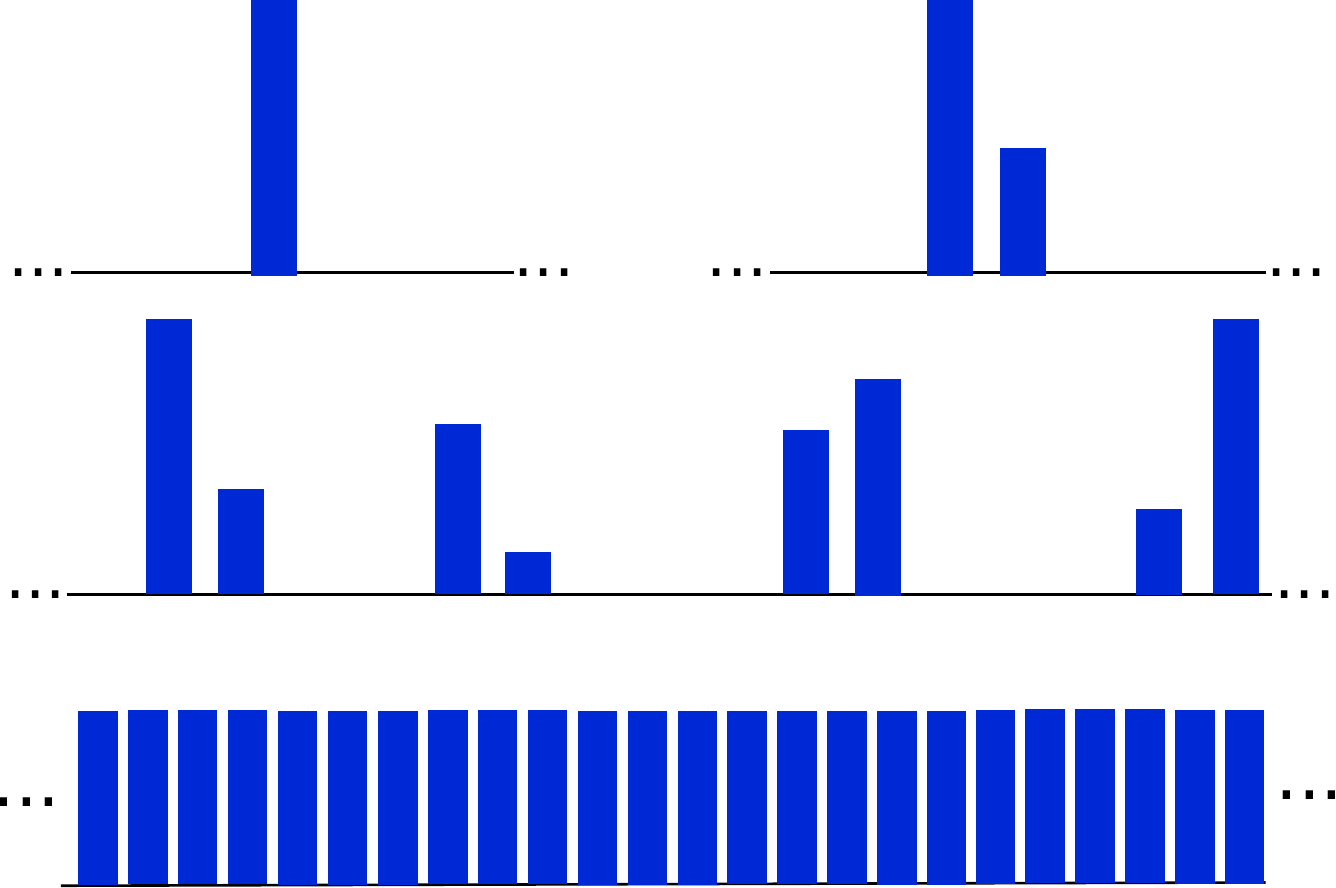}   \hspace{.5 in} \includegraphics[scale=0.4]{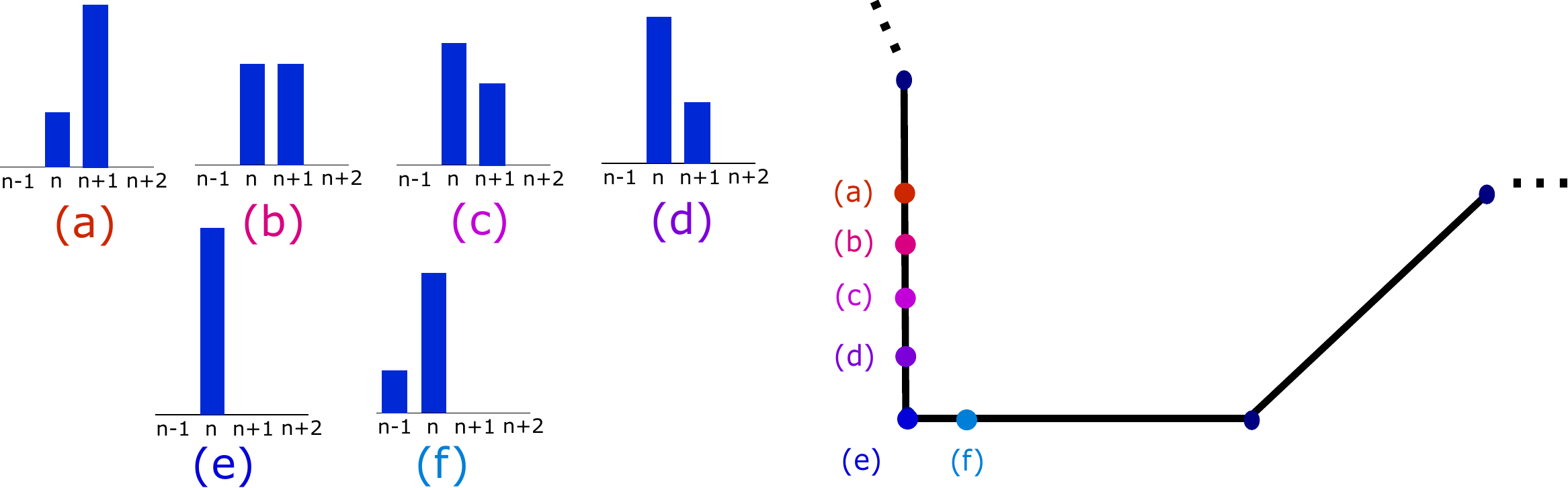}
    \caption{Single-party at one or two sites, multi-party, and uniform distribution equilibria (left). Phase-space schematic of a family of single-party equilibria with corners at one-site parties (right).}
    \label{fig:my_label}
\end{figure}

\subsection{Single party states}

A single two-site party of mass $m$ takes the form $P^*_{\alpha, n} = (0, ..., \overset{\downarrow n}{m\alpha}, m(1-\alpha), ... , 0)^\perp,$ for $ \alpha \in (0,1).$ In the limit $\alpha=0$ or $\alpha=1$, the party only occupies one site. It represents (total) consensus of opinion at two adjacent opinion sites $n$ and $n+1$. The linearized vector field at a two-site and at a one-site equilibrium ($\alpha=0$) of \eqref{e:bc} with mass $m=1$  are, respectively, given by the matrices
\[
\begin{pmatrix} 0 & \ &\ &\ &\ &\ &\ & \ &\ &\ \\
\ & \ddots & \  &\ &\ &\ &\ &\ &\  \ &\ \\ 
\ & \ & -\alpha &           0 & 0 & 0 & 0       & 0           & \ &\ \\
\ & \ & 2\alpha & -(1-\alpha) & 0 & 0 & 0       & 0           & \ &\ \\
\ & \ & -\alpha & 2(1-\alpha) & 0 & 0 & -\alpha & 0           & \ & \xleftarrow{n} \\
\ & \ &       0 & -(1-\alpha) & 0 & 0 & 2\alpha & -(1-\alpha) & \ &\ \\
\ & \ &       0 &           0 & 0 & 0 & - \alpha & 2(1-\alpha)& \ &\ \\
\ & \ & 0       &           0 & 0 & 0 & 0 & -(1-\alpha) & \ &\ \\
\ & \ & \   & \ &\ &\ &\ &\ &\ \ddots & \  \\
\ & \ & \ &\ & \underset{n}{\uparrow} \normalsize &\ &\ &\ & \ & 0 \end{pmatrix},\ \text{ and }
\ 
\begin{pmatrix} 0 & \ &\ &\ &\ &\ & \ &\ &\ \\
\ & \ddots & \  &\ &\ &\ &\ &\ &\  \  \\ 
\ & \ & -1 &           0 & 0 & 0      & 0           & \ &\ \\
\ & \ & 2 & 0 & 0 & 0 & 0                 & \ &\ \\
\ & \ & -1 & 0 & 0 & 0 & -1        & \ &\xleftarrow{n} \\
\ & \ &       0 & 0 & 0 & 0 & 2  & \ &\ \\
\ & \ &       0 &           0 & 0 & 0 & - 1 & \ &\ \\
\ & \ & \   & \ &\ &\ &\ &\  \ddots & \  \\
\ & \ & \ &\ & \underset{n}{\uparrow} \normalsize &\ &\ &\ & 0 \end{pmatrix}.
\]
By scaling invariance, linearization at parties with mass $m$ gives the same matrices multiplied by a factor $m$. 

Both matrices clearly have infinite-dimensional kernels, with bases spanned by 
$$ \left\{e_j \  \big| \  j \neq n-2, n -1, n+2, n+3\right\},\quad   \text{ and } \quad 
\left\{e_j \  \big| \  j \neq n -2, n+2\right\},
$$
respectively. Here, $e_n$ is the canonical basis vector such that $(e_n)_k= \delta_{nk}$ with Kronecker-$\delta$ notation. 

For the two-site party, the kernel has codimension 4, and corresponds to the two sites of the party and the sites at a distance 3 or more away. The spectrum of the linearization is $\{0,-m\alpha, m(1-\alpha)\}.$ For the one-site party, the kernel has  codimension 2 and the spectrum of the linearization is $\{0, -m \}$. The complement of the kernel can be associated with opinion sites $n$ that interact with the support of the party, while kernel elements correspond to lattice sites that do not interact with the party, either because they are too far away or because they simply change the shape of the party. 

Associated with the kernel elements $e_n$ and $e_{n+1}$ for the two-site party, there is a two-dimensional family of single party equilibria parameterized by the mass $m$ and the parameter $\alpha$, which can be thought of as parameterizing the average opinion in the party. At the two-site party, one can associate the kernel vector $e_n$ with mass change, and the kernel vectors $e_{n+1}$ and $e_{n-1}$ with increases or decrease in the average opinion in the party. Note however that the direction $e_{n+1}+e_{n-1}$  associated with states $(\ldots,0,0,\alpha,1,\alpha,0,0,\ldots)^\perp$ does not correspond to the tangent space of a family of equilibria. Similarly, directions in the kernel with support on more than two sites do not correspond to families of equilibria. In particular, the infinite-dimensional kernel is \emph{not} simply the tangent space to a high-dimensional family of equilibria, a fact that will slightly complicate the application of singular perturbation theory, later. We remark that the stability of single-party states is analytically rather subtle due to this high-dimensional kernel and the possible associated dynamics of clustering of small mass nearby in phase space, but far away in the opinion spectrum. 

Fixing total mass, single-party states naturally come in a one-parameter  family $\mathcal{F}$ that can be parameterized by their average opinion: A party with mass $1-\alpha$ at site $n$ and mass $\alpha$ at site $n+1$ has average opinion $n+\alpha$.

Drifting opinion in a biased model is to leading order described by drift along this continuous family of single-party states. A subtlety arises when viewing this family in phase space. The tangent vector to the family of two-site parties supported on sites $n$ and $n+1$ is $e_n-e_{n+1}$. This tangent vector is discontinuous at the one-site party, where the continuous curve of single-party equilibria possesses a corner; see Figure \ref{fig:my_label} for an illustration. Drift along this corner, as we shall see below, introduces dynamics and error terms known from the analysis of a  passage through a transcritical bifurcation.




\subsection{Uniform distribution of opinions}

\begin{figure}
\centering
\includegraphics[scale=0.75]{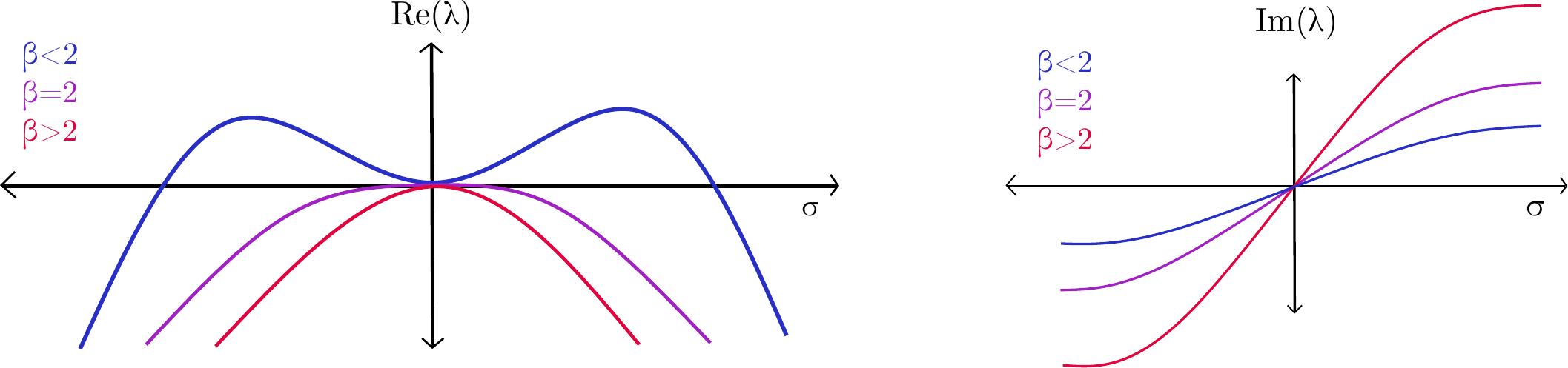}
\caption{Linear dispersion relation \eqref{e:dis} for $\beta<2$, $\beta=2$, and $\beta>2$.}\label{fig:dispersion}
\end{figure}

The uniform state $P_n \equiv m$ is also an equilibrium of \eqref{e:bc}. This equilibrium turns out to be unstable and a typical question of interest is how fluctuations around this equilibrium evolve into multi-party states. In order to understand this process, one usually starts by linearizing \eqref{e:bc} at the uniform state to find
\begin{equation*}
    \frac{dP_n}{dt} = m(-P_{n-2} + 2P_{n-1} -2P_n + 2P_{n+1} - P_{n+2}).
\end{equation*}
Solutions to this constant-coefficient lattice-differential equation can readily be found after Fourier transform. Therefore inserting an ansatz $P_n=e^{i\sigma n + \lambda t}$, one finds the dispersion relation
\begin{equation*}
    \lambda = 2m(2\cos(\sigma) - \cos(2\sigma) - 1), \qquad -\pi\leq \sigma<\pi;
\end{equation*} 
see also Figure \ref{fig:dispersion}.
The temporal eigenvalue $\lambda$ is real  and obtains a maximum of 1 at $\sigma = \pm \frac{\pi}{3}$. 

The linearization therefore predicts fastest growth of perturbations with period $n=6$, predicting that white-noise fluctuations around a constant state would evolve towards a multi-party state with party peaks at sites with distance $\delta n=6$. A more refined branch point analysis of this dispersion relation reveals that localized perturbations of the unstable state evolve into parties with different spacing, $\delta n=5.311086\ldots$; see \cite{bennaimscheel2016}. We will return to this analysis when considering stability of uniform states with (strong) bias. 


\subsection{The effect of bias on equilibria}

Returning to the model equation with self-incitement bias,
\begin{equation}\label{e:bc_qb} \frac{dP_n}{dt} = 2P_{n+1}P_{n-1} - P_n(P_{n+2} + P_{n-2}) + \beta(P^2_{n+1}-P_n^2),\end{equation} 
we note that, for $\beta>0$, single-party states do not form equilibria. The resulting drift of single parties along the family of equilibria $\mathcal{F}$ with a resulting change in average opinion is the object of much of the remainder of this paper. 

On the other hand, the uniform state does persist as an equilibrium, but the linearization picks up new terms. When $\beta >0$, the linearized equation
\begin{equation*} \frac{dP_n}{dt} = m(-P_{n-2} + 2P_{n-1} -(2+2\beta)P_n + (2+2\beta)P_{n+1} - P_{n+2})
\end{equation*}
has dispersion relation
\begin{equation}\label{e:dis} \lambda = m(-2\cos(2\sigma) + (4+2\beta)\cos(\sigma) -(2+2\beta) + i\beta\sin(\sigma)).\end{equation}
For $\beta < 2$, the maximum of $\Re\lambda(\sigma)$ is positive  and the uniform state is unstable. However, for  $\beta \ge 2$,  $\Re\lambda\leq 0$ is nonpositive, with a quadratic tangency of the eigenvalues at the origin for $\beta>2$: strong bias stabilizes uniform distribution of opinion and hence disfavors consensus! 

The spectral stability in the dispersion relation translates readily into linear stability in, say, $\ell^2$. One would also expect nonlinear stability with approximately diffusive decay of the perturbation for localized initial conditions, measured in $\ell^\infty$, due to the presence of (discrete) derivatives in the nonlinearity; see Remark \ref{r:bst} for more detail. 

We shall exploit the change of stability in our bifurcation analysis, showing that the destabilization is accompanied by the creation of localized coherent structures, in \S \ref{s:large}.

\section{Small bias regime}\label{s:small}

We study here dynamics for $0<\beta \ll 1$. We use methods from geometric singular perturbation theory (GSPT) to investigate the speed of propagation of a single party in the biased system. We note that the lack of smoothness of the family of equilibria in the unbiased system prohibits a global slow-fast decomposition using existing theory. Nevertheless, we find locally invariant manifolds and separately analyze the system near the points where the manifold is not smooth, allowing us to compute the local speed of propagation for small bias. 

\subsection{Singular perturbation analysis near single-party equilibria}

We  begin by establishing the existence of locally invariant manifolds. 

\begin{Proposition}[Two-site center manifold]\label{p:2site}
    Fix $n \in \Z$, $k\geq 1 \in \N$ arbitrary, and $\delta > 0$, arbitrarily small. Then, for each fixed mass $m>0$ of two-site party equilibria,  there exists a family of locally invariant, infinite-dimensional, codimension-4,  $C^k$-manifolds $\mathcal{M}_\beta\subset\ell^\infty$, which depend on $\beta$ in a  $C^k$-fashion  such that $\mathcal{M}_0 \supset \{P^*_{\alpha,n} \ | \ \alpha \in (\delta, 1 - \delta) \}$, the part of the family of two-site single party equilibria supported on sites $n$ and $n+1$ away from one-site party equilibria. Its tangent space for $\beta=0$ at any of the two-site equilibria coincides with the infinite-dimensional kernel of the linearization at this equilibrium. 
\end{Proposition}




\begin{Proof}
The manifold of equilibria $\mathcal{M}_0$ is (locally) invariant and its linearization possesses an exponential dichotomy with a 4-dimensional stable subspace and an infinite-dimensional center subspace, with uniformly bounded projections $\mathcal{P}^s$ and $\mathcal{P}^c$. Standard theory for invariant manifold then shows the existence of smooth, locally invariant center-manifolds associated with this splitting, using for instance graph transforms as in \cite{batesluspike2008, fenichel1,hirschpughshub}.
\end{Proof}

\begin{Proposition}[One-site center manifold]\label{p:1site}
   Fix $n \in \Z$, $k\geq 1 \in \N$ arbitrary, and $\delta > 0$, sufficiently small. Then, for fixed mass $m>0$ of the one-site equilibrium,  there exists a family of locally invariant, infinite-dimensional, codimension-2,  $C^k$-manifolds $\mathcal{M}'_\beta\subset\ell^\infty$, which depend on $\beta$ in a  $C^k$-fashion  such that $\mathcal{M}'_0 \supset  P_n^* \cup \{ P^*_{\alpha,m} \ | m = n, \alpha \in (0, \delta);  m = n-1, \alpha \in (1-\delta,1) \}$, the family of two-site single party equilibria close to the one-site equilibrium $P^*_n$. Its tangent space for $\beta=0$ at the one-site equilibrium coincides with the infinite-dimensional kernel of the linearization at this equilibrium. 
\end{Proposition}

\begin{Proof}
    This is a standard local center-manifold result in infinite dimensions; see for instance \cite{VI,Henry1989GeometricTO}. It contains the family of equilibria since it contains all small solutions bounded for all times. 
\end{Proof}

We emphasize that the two-site center-manifold is global in the sense that it contains a compact subset of the line of equilibria between two one-site parties, while the one-site center manifold is local, defined only in a small neighborhood of the one-site party. 

\subsection{Leading-order dynamics away from corners}
Using invariance, we can now compute the leading-order dynamics on each manifold and the reduced flow. Away from the corners, 
we can parameterize the kernel of the linearization at a two-site party $P^*_{n,\gamma}$ through values $q_n$, $q_{n+1}$, and $p\in \ell^\infty,$ $p_j=0$ for $n-2 \le j \le n+3$, and write the manifold by $P=h(q_n,q_{n+1},p,\beta)$, with 
\[
h(P_n,P_{n+1},p,\beta)= (\gamma + q_n)e_n+(1 - \gamma + q_{n+1})e_{n+1}+p+\mathcal{O}(\beta,|q_n|^2+|q_{n+1}|^2,|p|^2).
\]
In this parameterization, we have
\begin{align*}
\dot{q}_n&=\beta \mathcal{P}^c(G(h(q_n,q_{n+1},p,\beta)))_n + \mathcal{O}(\beta^2,|q_n|^2+|q_{n+1}|^2,|p|^2),\\
\dot{q}_{n+1}&=\beta \mathcal{P}^c(G(h(q_n,q_{n+1},p,\beta)))_{n+1} + \mathcal{O}(\beta^2,|q_n|^2+|q_{n+1}|^2,|p|^2),\\
\dot{p}&=\mathcal{O}(|p|^2,|\beta p|,\beta^2(|q_n|+|q_{n+1}|)),
\end{align*} 
where $G(P)_j = P_{j+1}^2 - P_j^2,$ and $\mathcal{P}^c$ is the spectral projection onto $\ker Df(P^*_{\alpha,n})$, defined as

\[ \mathcal{P}^c(P)_j = \begin{cases}  P_j, & j < n-2, j > n+3 \\ 0 , & j = n+1,n+2, n-1,n-2 \\ 3P_{n-2} + 2P_{n-1} + P_n  - P_{n+2} - 2P_{n+3}, & j = n\\ -2P_{n-2} - P_{n-1} + P_{n+1} +2 P_{n+2} + 3 P_{n+3}, & j = 
n+1\end{cases}, \]
Now scaling $p=\beta q$, we find explicitly at leading order on the center manifold, 
\begin{align*}
\dot q_n&=  \beta(\gamma^2 + (1-\gamma)^2)+\mathcal{O}(\beta^2),\\
\dot q_{n+1}&=  -\beta((\gamma)^2 + (1-\gamma)^2)+\mathcal{O}(\beta^2),\\
\dot{p}&= \mathcal{O}(\beta^2).
\end{align*}
Clearly, at leading order mass in the party is conserved, $\dot P_n+\dot P_{n+1}=0$. For the family of two-site parties parameterized by $P_n=\alpha, P_{n+1}=1-\alpha$, this gives
\begin{equation}\label{e:drift}
\dot{\alpha}=-\beta[\alpha^2+(1-\alpha)^2]+\mathcal{O}(\beta^2).
\end{equation}

\begin{figure}[H]
    \centering
     \includegraphics[width=.6\textwidth]{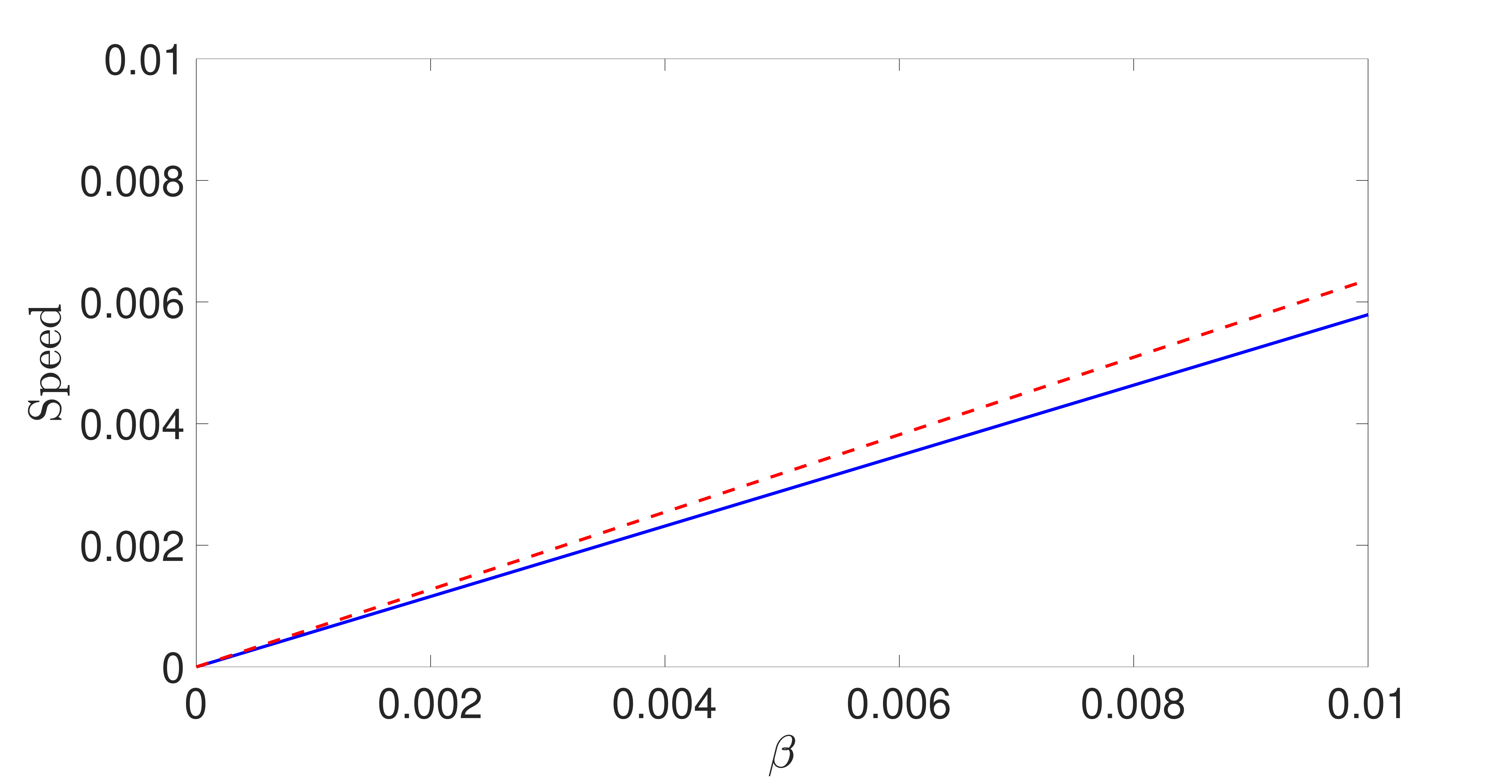}
    \caption{Comparisons of leading order average drift speed prediction from theory (red) with numerically computed drift speed (blue) for a party of total mass 1.}
    \label{fig:gstp-speeds}
\end{figure}

We compare this first order approximation \eqref{e:drift} with numerically computed drift speeds averaged in time, that is, covering the interval $\alpha\in [0,1]$ rather than the interval $(\delta,(1-\delta)$ where the above analysis applies, in Figure \ref{fig:gstp-speeds}. Details on computational procedures are  delineated in \S{}\ref{s:num}. We see that the predicted speeds here give the leading order term for small $\beta$, but the discrepancy for even moderately small values of $\beta$ are significant. We show how this discrepancy can be attributed to the passage near the one-site parties, contributing a term $\beta^{3/2}$. 

\subsection{Corner dynamics and slow passage through transcritical bifurcations}
We now compute the leading-order dynamics near the one-site party. The kernel of the linearization $Df(P_n^*)$ is now of codimension 2, and we parameterize elements $P_c$ of the kernel by values $\alpha_+, \alpha_-, \alpha_m \in \R$ and $p \in \ell^\infty$, with $p_j = 0$ for $n-2 \le j \le n+2$. Let 
\begin{equation*}
\displaystyle
 \begin{array}{lll} 
     e_+ = (0, ... , 1, \overset{\downarrow n}{-1}, 0, ...)^\perp,& e_- = (0, ..., 0, \overset{\downarrow n}{-1}, 1, ...)^\perp, &   e_m = (0, ... , 0, \overset{\downarrow n}{1}, 0, ...)^\perp,\\
        e_+^* = (0, ... , 0,0, \overset{\downarrow n}{0}, 1, 2 ...)^\perp, &e_-^* = (0, ..., 2,1,\overset{\downarrow n}{0}, 0, 0, ...)^\perp,& e_m^* = (0, ..., 1,1,\overset{\downarrow n}{1}, 1, 1, 0 ...)^\perp.
\end{array}
\end{equation*}

We write the  center manifold as manifold as the graph of $h^c= h^c(\alpha_+, \alpha_-, \alpha_m, p, \beta)$, with 
\[
h^c(\alpha_+, \alpha_-, \alpha_m, p, \beta) = 
 P_n^* + p + \alpha_+e_+ + \alpha_-e_- + \alpha_m e_n  + \mathcal{O}(| \alpha_+|^2 + |\alpha_-|^2 + |\alpha_m|^2 + |p|^2 + \beta^2).
\]
 In these coordinates, we have 
\[ \dot{P}_c = \mathcal{P}^c(f(P_n^* P_c + h^c(P_c))) + \mathcal{O}(|P_c|^3 + \beta^3),\]
where the spectral projection $\mathcal{P}^c$ onto $\ker Df(P^*_{\alpha,n})$ is given by \[ \mathcal{P}^c v = \langle v, e_+^* \rangle e_+ + \langle v, e_-^* \rangle e_- + \langle v, e_m^* \rangle e_n + \sum_{|i|>3} \langle v, e_i \rangle e_i.\]

Writing the vector field in terms of $\alpha_+, \alpha_-, \alpha_m, p, \beta$, we find
\[f(P_n^* + P_c + h(P_c)) = \begin{pmatrix}
\vdots  \\
2p_{-3}p_{-5} - p_{-4}p_{-6}\\ 
-p_{-3}\alpha_- - p_{-3}p_{-5}\\ 
-h_- + 2p_{-3}\alpha_-\\
2h_- \alpha_-\alpha_+ - \alpha_-p_{-3} + \beta(1+2(\alpha_m - \alpha_- - \alpha_+))\\
2\alpha_+\alpha_- - (h_- + h_+) - \beta(1+2(\alpha_m - \alpha_- - \alpha_+)) \\
2h_+ - \alpha_+\alpha_- - \alpha_+p_3\\
-h_+ + 2\alpha_+p_3\\
-\alpha_+p_3 - p_3p_5 \\
2p_3p_5 - p_4p_6 \\
\vdots  \\ 
\end{pmatrix}\xleftarrow{n}  \qquad + \mathcal{O}(|P_c + \beta|^3),\]
which in turn gives
\begin{align*}
\dot{\alpha}_- &= \langle f(P_n^* + P_c + h^c(P_c)), e_-^* \rangle &&= \ \ \ \beta - \alpha_-\alpha_+ + 3p_{-3}\alpha_- + 2\beta(\alpha_m - \alpha_+ - \alpha_-) + \mathcal{O}(3)\\
\dot{\alpha}_+ &= \langle f(P_n^* + P_c + h^c(P_c)), e_+^* \rangle &&= \ -\alpha_-\alpha_+ + 3 p_3\alpha_+ + \mathcal{O}(3)\\
\dot{\alpha}_m &= \langle f(P_n^* + P_c + h^c(P_c)), e_m^* \rangle &&= \ \  \ p_3\alpha_+ + p_{-3}\alpha_- + \mathcal{O}(3)\\
\dot{p}_3 &= \langle f(P_n^* + P_c + h^c(P_c)), e_3 \rangle &&= \ -p_3\alpha_+ - p_3p_5 + \mathcal{O}(3)\\
\dot{p}_{-3} &= \left\langle f(P_n^* + P_c + h^c(P_c)), e_{-3} \right\rangle &&= \  -p_{-3}\alpha_- - p_{-3}p_{-5} + \mathcal{O}(3)\\
\dot{p}_{4} &= \left\langle f(P_n^* + P_c + h^c(P_c)), e_{4} \right\rangle &&= \  - 2p_{3}p_{5} - p_{4}p_{6} + \mathcal{O}(3)\\
\dot{p}_{-4} &= \left\langle f(P_n^* + P_c + h^c(P_c)), e_{-4} \right\rangle &&= \  - 2p_{-3}p_{-5} - p_{-4}p_{-6} + \mathcal{O}(3) \\
\dot{p}_j &= 2p_{j-1}p_{j+1} - p_j(p_{j+2}p_{j-2} + \mathcal{O}(3), &&|j| > 4.
\end{align*}
At leading order, the subspace where $p_j=0$ for all $j$ and $\alpha_m=0$ is invariant, and we therefore consider the leading-order equation for $\alpha_-,\alpha_+$, only,
\[
\dot{\alpha}_- = \beta - \alpha_-\alpha_+  - 2\beta( \alpha_+ + \alpha_-),\qquad \dot{\alpha}_+ =  -\alpha_-\alpha_+ .
\]
Changing variables  $x = \alpha_- + \alpha_+$ and $\mu = \alpha_- - \alpha_+$, this gives
\[
\dot{x}=-\frac{1}{2}(x^2-\mu^2)+2\beta(1-x),\qquad \dot{\mu}=2\beta(1-x).
\]
In the natural scaling $x\sim\mu\sim \beta^{1/2}$, the terms $-2\beta x$ are of higher order. The remaining terms describe precisely the slow passage through a transcritical bifurcation as studied for instance in \cite{krupa2001transcritical} using geometric desingularization. 

\begin{figure}[H]
    \centering
     \includegraphics[width=.5\textwidth]{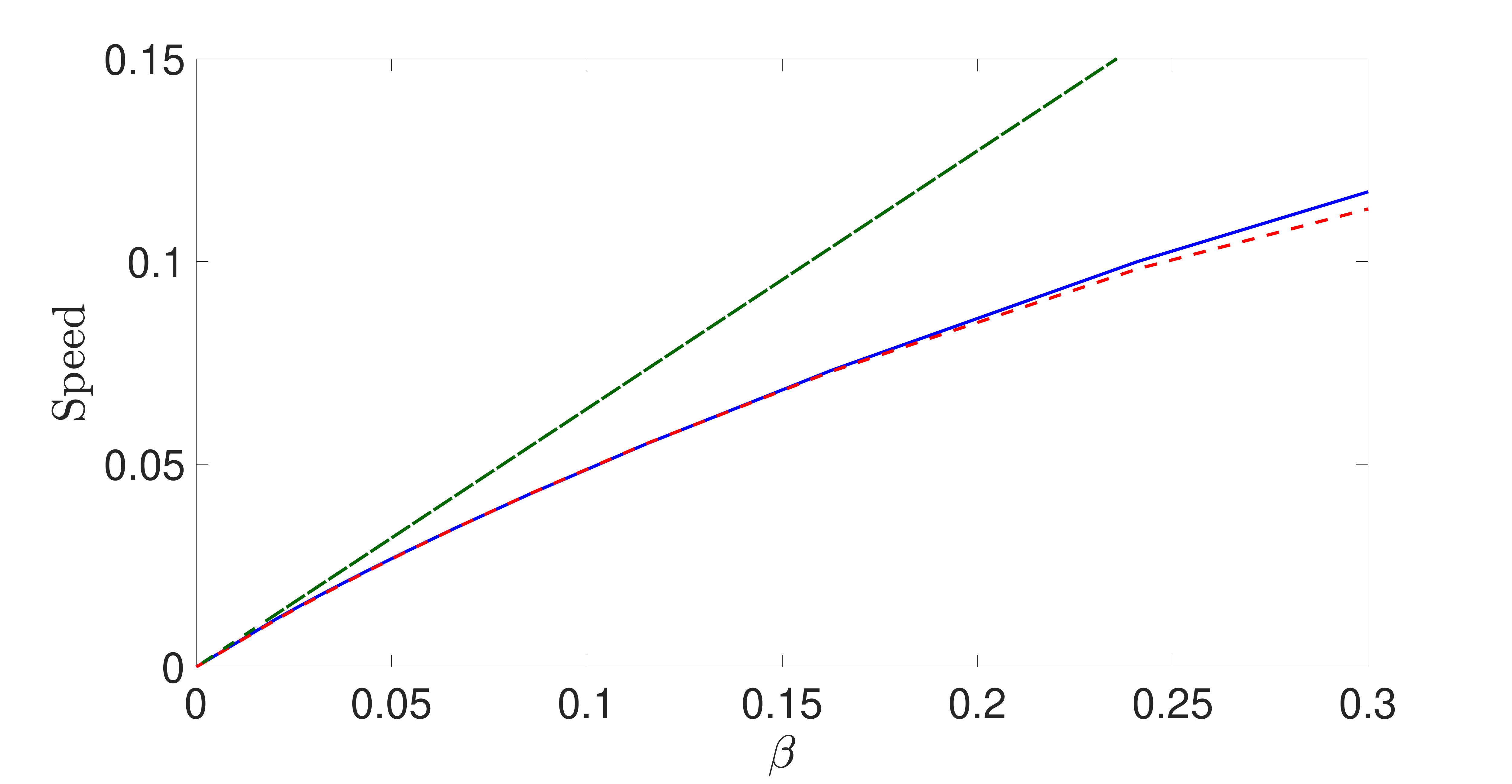}
    \caption{Numerically computed values for the drift speed (blue) plotted against the leading-order approximation (green) and the approximation with the $\beta^\frac{3}{2}$ correction.}
    \label{fig:threehalf}
\end{figure}

The passage near the one-site party can be analyzed by starting in a section to the flow $\{x-\mu=2\delta\}$ near $x=-\mu=\delta$ and tracking time until the section  $x+\mu=2\delta$ near $x=\mu=\delta$. At leading order this time is given by $2\delta/\beta$, which confirms that, at leading order, the passage time near the one-site equilibrium can be ignored in the computation of the average speed, letting for instance $\delta\to0$. The scaling does however introduce error terms involving $\beta^{1/2}$, which indeed manifest themselves in corrections to the averaged speed of order $\beta^{3/2}$; see for instance \cite{krupa2001transcritical}. We did not attempt to derive these corrections analytically but 
numerically found the coefficient to the $\beta^{3/2}$ correction  as -.467; see  Figure \ref{fig:threehalf} for numerical values of the drift speed of a party of mass 1, plotted against the theoretical prediction with and without the $\beta^{3/2}$-correction. 


\section{Large bias regime}\label{s:large}

In this section, we investigate the dynamics as $\beta$ gets large. Numerically, we see that the background mass of profiles increases, and that as $\beta$ approaches 2, the size of the profiles becomes arbitrarily small relative to the background mass. 

We can see that this agrees with heuristics from the spectrum of the linearization at the background mass. We recall the dispersion relation for the linearization at the uniform steady-state $u \equiv m$:
\[ -i\omega = m(-2\cos(2\sigma) + (4+2\beta)\cos(\sigma) -(2+2\beta) + 2i\beta\sin(\sigma)),\] 
and note that the quantity $-i\omega$ has nonpositive real part exactly when $\beta \ge 2$. Therefore the uniform steady-state regains (marginal) linear stability when $\beta = 2$, which may explain why traveling profiles disappear. For $\beta < 2$ but sufficiently close, numerically we see small-amplitude, long-wavelength profiles atop the background state $u \equiv m$. 
In this section, we rigorously establish the existence of these spike solutions for $\beta$ sufficiently close to 2. We first derive formal amplitude equations, giving heuristics for the existence of spike solutions. We then establish  existence rigorously  and derive expansions through a nonlocal center manifold computation, followed by a Melnikov analysis of the reduced system.

\paragraph{Formal derivation of amplitude equations.}
We derive formal amplitude equations for the $\beta \sim 2$ regime, under long-wavelength, small-amplitude assumptions. We will see in \S{}\ref{ss:cm} that rigorous center manifold calculations in fact agree with the resulting amplitude equations. For our approximation, we follow \cite{sw00fpu} by choosing the regime in which the lattice spacing remains constant, but the spatial variable in the KdV equation is rescaled. 
\par Assume that \eqref{e:bc_qb} permits a small, long-wavelength solution $P_n = m + \ep^2 u(\ep(n+ct), \ep^3 t)$ where $c = 2 m \beta$ is chosen as the group velocity at the constant state. To simplify notation, let $\xi \equiv n+ct$. Substituting the ansatz into \eqref{e:bc_qb} and grouping linear and nonlinear terms, we find
\begin{align}\label{ueqLong} \partial_t \big[u (\ep\xi, \ep^3 t) \big] = &2m \big[  u ( \ep\xi - \ep, \ep^3 t) + u (\ep\xi + \ep, \ep^3 t) \big] - m\big[2  u (\ep\xi, \ep^3 t) +  u (\ep\xi -2\ep, \ep^3 t) +  u (\ep\xi+2\ep, \ep^3 t) \big]\nonumber \\
&+ 2m \beta \left[ u (\ep\xi + \ep, \ep^3 t)  -   u (\ep\xi , \ep^3 t) \right] -\ep^2 u (\ep\xi, \ep^3 t) \big[u (\ep\xi -2\ep, \ep^3 t) +  u (\ep\xi+2\ep, \ep^3 t) \big] \nonumber\\
& + 2 \ep^2  u ( \ep\xi - \ep, \ep^3 t)  u (\ep\xi + \ep, \ep^3 t) + \beta \ep^2 \left[ u^2 (\ep\xi + \ep, \ep^3 t)  -   u ^2(\ep\xi , \ep^3 t) \right].\end{align}
We now expand $u$ in the first variable,
$$ u(\ep\xi + k \ep, \ep^3 t) = u(\ep\xi, \ep^3 t) + k \ep \cdot \partial_1 u(\ep\xi, \ep^3 t)+ \frac{(k\ep)^2}{2!}\cdot \partial_1^2 u(\ep \xi, \ep^3 t) +  \frac{(k\ep)^3}{3!}\cdot \partial_1^3 u(\ep \xi, \ep^3 t) + \frac{(k\ep)^4}{4!}\cdot \partial_1^4 u(\ep \xi, \ep^3 t) + \mathcal{O}(\ep^5),$$
for each $k \in \{-2,-1,0,1,2\}$, and we implicitly assume sufficient smoothness in $u$. Inserting the expansion into \eqref{ueqLong}, we find 
\begin{align*} \ep^3 \partial_2 u  = &2m \big[\ep^2 \partial_1^2 u  + \frac{\ep^4}{12} \partial_1^4 u \big] - m\big[4  \ep^2 \partial_1^2 u  + \frac{4}{3} \ep^4 \partial_1^4 u  \big] + 2m \beta \left[ \frac{\ep^2}{2} \partial_1^2 u  + \frac{\ep^3}{6} \partial_1^3u  + \frac{\ep^4}{24} \partial_1^4 u  \right] \\
&-\ep^2  \big[2u^2 + 4 \ep^2 u  \partial_1^2u  + \frac{4 \ep^4}{3} u \partial_1^4u  \big]  + 2 \ep^2  \Big[ u^2  + \ep^2 u  \partial_1^2u  - \ep^2 (\partial_1 u)^2 + \frac{\ep^4}{12} u  \partial_1^4u - \frac{\ep^4}{3} \partial_1 u\partial_1^3u + \frac{\ep^4}{4} (\partial_1^2u)^2 \Big] \\
&+ \beta \ep^2 \Big[ 2 \ep u \partial_1 u  + \ep^2 u \partial_1^2u + \ep^2  (\partial_1 u)^2  + \ep^3 \partial_1 u \partial_1^2 u + \frac{\ep^3}{3} u \partial_1^3 u + \frac{\ep^4}{4} (\partial_1^2 u)^2  + \frac{\ep^4}{12} u \partial_1^4 u + \frac{\ep^4}{3} \partial_1 u \partial_1^3u  \Big] + \mathcal{O}(\ep^5) .\end{align*}
where $\partial_1$ and $\partial_2$ refer to the partial derivatives with respect to the first and second arguments $\ep \xi$ and $\ep^3 t$, respectively. Retaining orders $\mathcal{O}(\ep^2)$ to $\mathcal{O}(\ep^4)$, then dividing by $\ep^3$, we find the formal amplitude equation 
\begin{align}\label{eq:mod}
 \partial_2 u = \frac{m(\beta-2)}{\ep} \partial_1^2 u  + \frac{m \beta}{3} \partial_1^3 u + 2\beta u \partial_1 u- \ep \left( \frac{(14-\beta)m}{12} \partial_1^4 u +(2-\beta) \partial_1 \left( u \partial_1 u\right)  \right) + \mathcal{O}(\ep^2).
\end{align}
Note that for $\beta = 2 + \mathcal{O}(\ep^2)$, at leading order, \eqref{eq:mod} becomes the KdV equation
\begin{equation}\label{eq:KdV1}
    \partial_2 u  =    \frac{2m}{3} \partial_1^3 u + 4 u \partial_1 u .
\end{equation}
The KdV equation possesses a family of traveling spikes parameterized by the wave speed (or the amplitude). 
At the next order $\mathcal{O}(\ep)$, with $\beta=2-\ep^2\tilde{\beta}$, $\tilde{\beta}>0$, we find a damping term $-\ep  m \partial_1^4 u$, fourth order viscosity, and a negative damping term $-\ep m\tilde{\beta}\partial_1^2 u$, negative viscosity. The analysis presented below demonstrates that this equation can be rigorously derived as an ODE at leading order on a center manifold, and that the effects of negative and positive viscosity balance for an appropriate wave speed (or amplitude) in the reduced center manifold equation.

\begin{Remark}[Viscous Burgers modulation and stability]\label{r:bst}
For $\beta>2$ one finds that \eqref{eq:mod} reduces at leading order to viscous Burgers equation, $\partial_2 u = \partial_1^2 u + u \partial_1 u$, after appropriate scalings. It is  known that in this approximation, localized initial conditions decay algebraically in $L^\infty$  and that higher-order terms, that is, terms carrying higher powers in $u$ or more derivatives are irrelevant in the long-time asymptotics of small data; see \cite{bk94}. We therefore suspect that one can establish asymptotic stability in  $\ell^\infty$ of constant distributions of opinions in our system for values of $\beta>2$ and small perturbations in $\ell^1$. On the other hand, $\ell^\infty$ perturbations may evolve into persistent dynamics, for instance viscous shocks, rarefaction waves, and their superposition; see \cite{dsss} for such a construction based on a Burgers approximation. 
\end{Remark}



\subsection{Statement of main result} In order to precisely state the main result of this section, we collect some definitions and notation. We define a traveling wave in a 1-D lattice as a solution of the form $P_n = Q(n+ct)$, where $Q(\cdot)$ defines a fixed profile which moves through the lattice. Seeking such a solution, we reformulate the problem as a functional differential equation: substituting $P_n=Q(n+ct)$, we get
\begin{equation}\label{e:FDE}
   0 = - c Q'(\xi) + 2 Q(\xi-1) Q(\xi+1) - Q(\xi) \big( Q(\xi-2) + Q(\xi+2) \big) + \beta \big( Q^2(\xi+1) - Q^2(\xi) \big),
\end{equation}
where $\xi=n+ct \in \R$.

Further writing $Q(\xi) = m + q(\xi)$, we have
\begin{equation}\label{e:fdefull}
        0 =  -c q'(\xi) + m(2q(\xi -1) + 2(\beta + 1)q(\xi + 1) - 2(\beta + 1)q(\xi) - q(\xi -2) - q(\xi+2) + \mathcal{N}(q,\beta)),
\end{equation}
where $\mathcal{N}(q,\beta) = 2 q(\xi-1)q(\xi+1) - q(\xi) ( q(\xi-2) + q(\xi+2)) + \beta ( q(\xi+1)^2 -  q(\xi)^2)$. 

Note that we can rewrite the forward-backward delay equation \eqref{e:fdefull} as a nonlocal equation
\begin{equation}\label{e:FDEq}
    0 = -c q'(\xi) - 2m (\beta + 1) q(\xi) + m \mathcal{K}*q + \mathcal{N}(q,\beta),
\end{equation}
with convolution kernel $\mathcal{K}(\cdot) = -\delta(\cdot-2) + 2 \delta(\cdot-1) + 2(\beta+1) \delta(\cdot+1) - \delta(\cdot+2) $.

We are now ready to state the main result of this section:

\begin{Theorem}\label{existence_spikes}
There exists $\beta_* > 0$ such that for $2-\beta_* < \beta < 2$, there exists a locally unique homoclinic solution $q=q^*_\beta(\xi)$ to \eqref{e:FDEq} with locally unique wave speed $c=c(\beta)$, with $\lim_{|\xi| \to \infty}q^*_\beta(\xi) = m$, and therefore a unique traveling-wave spike solution $({u_\beta^*})_n(t) = q_\beta^*(n +ct)$ to \eqref{e:bc_qb}. Furthermore, we have \begin{equation*} \begin{split}
    q_\beta^*(\xi) &= m( 1 + \frac{7(2-\beta)}{10})\sech^2(\sqrt{\frac{7(2-\beta)}{20}}\xi) + \mathcal{O}((2-\beta)^\frac{3}{2})\\
    c(\beta) & = m(4 - \frac{16}{15}(2-\beta)) + \mathcal{O}((2-\beta)^\frac{3}{2}).
\end{split}
\end{equation*}
Note that in particular $q_\beta^*(\xi)>m$ represents an opinion spike, a traveling party, and $c=c(\beta)>2m\beta$ is slightly larger than the linear group velocity $2m\beta$.
\end{Theorem}

The outline of the proof is as follows: We first show that \eqref{e:FDEq} can be reformulated to satisfy the hypotheses of the nonlocal center manifold theorem in \cite{fs2018center}, allowing for a center manifold reduction. We then calculate the nonlocal center manifold expansion in function space and derive the reduced vector field. Finally, we prove existence of homoclinic solutions to the reduced equations using Melnikov analysis, where the existence of a conserved quantity becomes crucial.

We begin with the existence of a center manifold. Let $\tbeta = 2 - \beta$, $\tc = c - 4m$, and define \begin{equation*}\T_{\tbeta, \tc}q = -cq' -2m(3 - \tbeta)q + m\mathcal{K}*q, 
\end{equation*}
as an operator on $H^1_{-\eta}(\R,\R)$ where 
\[
\|u\|_{H^1_{-\eta}}=\|u(\cdot) e^{-\eta |\cdot|}\|_{H^1},
\]
is an exponentially weighted space allowing for exponential growth. For $\eta>0$ sufficiently small, consider the kernel $\mathcal{E}_0$, which turns out to be finite-dimensional, and choose a closed complement, thus defining a projection $\mathcal{P}_0$ onto the kernel $\mathcal{P}_0$,
\begin{equation*}\mathcal{E}_0 = \ker\T_{0,0}, \qquad \mathcal{P}_0 \mathcal{E}_0=\mathcal{E}_0.\end{equation*}

\begin{Proposition}\label{p:cmfd}
Fix $\eta>0$ sufficiently small and consider the functional equation \eqref{e:FDEq} in $H^1_{-\eta}(\R,\R)$ and $k<\infty$. There exists neighborhoods $\mathcal{U}_q \times \mathcal{U}_{(0,0)}$ of $(0,0,0)$ in $\mathcal{E}_0 \times \R^2$ and a map $\Psi \in \mathcal{C}^k (\mathcal{U}_q \times \mathcal{U}_{(0,0)}, \ker \mathcal{P}_0 )$, with $\Psi(0,0,0) = 0$ and $D_q \Psi (0,0,0) = 0$ such that for all $(\tbeta, \tc) \in \mathcal{U}_{(0,0)}$ the manifold
$$ \mathcal{M}_0^{\tbeta,\tc} = \left\{ q_0 + \Psi(q_0, \tbeta, \tc): q_0 \in \mathcal{U}_q, \right\}\subset H^1_{-\eta} $$
contains the set of all bounded solutions of \eqref{e:FDEq}, small in $C^0(\R,\R)$.
\end{Proposition}

\begin{Proof} Preconditioning \eqref{e:FDEq} with the operator $\left( 2m (\beta+1) + c \frac{d}{d \xi} \right)^{-1} = G_{\beta,c}*$, with \[G_{\beta,c}(\xi) =  \frac{1}{c} \exp \left(- \frac{2m(\beta+1)}{c} \xi \right) \bigchi_{\mathbb{R}^+}(\xi),\] the resulting equation 
\begin{equation} \begin{split}
0 =& -q + (G_{\beta, c} * m \mathcal{K}) * q + G_{\beta, c}*\mathcal{N}(q,\beta)\\  =:& -q + \wt{\mathcal{K}}_{\beta,c}*q + \wt{\mathcal{N}}(q,\beta)\end{split}\end{equation}
is in the form of \cite{fs2018center}, with one discrepancy---the kernel $\wt{\mathcal{K}}_{\beta,c} = G_{\beta, c} * m \mathcal{K}$ is not in $ W^{1,1}_{\eta_0}$ for any $\eta_0 > 0$. However, since $\wt{\mathcal{K}}_{\beta, c}'$ is still a sum of an $L^1_{\eta_0}$ function and scaled translates of Dirac deltas, the results of \cite{fs2018center} can be established with minor modifications as follows. In the proof of \cite[Lemma 3.1]{fs2018center}, the smoothness of $\mathcal{K}$ is only used to show that $\mathcal{D}\T(u) = \mathcal{D}(\mathcal{K}*u) = (\mathcal{K}' + \rho \mathcal{K} + \rho \delta_0)*u $ satisfies the hypotheses of \cite{faye2014fredholm}. Here, we note that although  $\wt{\mathcal{K}}_{\beta,c} \notin  W^{1,1}_{\eta_0}$, it is still true that $(\wt{\mathcal{K}}_{\beta,c}' + \rho \wt{\mathcal{K}}_{\beta,c} + \rho \delta_0)*u$ satisfies the hypotheses of \cite{faye2014fredholm}. Therefore \cite[Lemma 3.1]{fs2018center} holds identically, and the rest of the proof in \cite{fs2018center} does not use this hypothesis. The parameter-dependent center manifold theorem \cite[Theorem 3]{fs2018center} can therefore be applied to the system \eqref{e:FDEq}, which implies the statement of the proposition. \end{Proof}

\subsection{Nonlocal center manifold expansion}\label{ss:cm}

Given the existence of a nonlocal center manifold, we now use the methods of \cite{fs2018center} to calculate the Taylor expansion of the center manifold in function space, and derive the reduced vector field. 

Since the Taylor expansion is written as a map over the kernel of $\T_{\tbeta,\tc}$, we first find a parametrization of $\ker\T_{\tbeta,\tc}$, and a projection. The dispersion relation given by the linearized equation, in terms of $\tbeta, \tc$, is
$$ d(\nu,\tbeta, \tc) = \widehat{\mathcal{T}_{\tbeta,\tc}}(\nu) = -(4m + \tc) \nu - 2m(3 + \tbeta) + 2m (e^{\nu} + e^{-\nu}) - 2m (e^{2\nu}+e^{-2\nu})+ 2m(2 - \tbeta) e^{\nu} = 0.$$
Note that $d$ is clearly analytic and roots on the imaginary axis $\nu\in i\R$ are a priori bounded. At $\tbeta=\tc=0$, we find that there are no imaginary roots $\nu\neq 0$.   
Expanding $d$ at $\nu =0$, we have 
$$ d(\nu,\tbeta, \tc) = (2m \tbeta -\tc) \nu + m\tbeta \nu^2 + \frac{m (2-\tbeta)}{3} \nu ^3 - \frac{m(\tbeta+12)}{12}\nu^4 + \mathcal{O}(\nu^5).  $$
Note that $\wt{d}(0,0,0) = \partial_{\nu} \wt{d}(0,0,0) = \partial_{\nu\nu}\wt{d}(0,0,0) = 0$, with $\partial_{\nu\nu\nu} \wt{d}(0,0,0) \neq 0$. Additionally, $\wt{d}(i\ell, 0,0) \neq 0$ for all $\ell \neq 0$. Thus the kernel $\mathcal{E}_0$ of $\mathcal{T}_{0,0}$ in $H^1_{-\eta}$ is given by
$$\mathcal{E}_0= \text{span}\{1, \xi, \xi^2 \},$$
and we write  elements $q_0 \in \mathcal{E}_0$ as
\begin{equation}\label{e:q_0} q_0(\xi) = A_0 + A_1 \xi + A_2 \xi^2 \in \mathcal{E}_0,\end{equation}
with $(A_0,A_1,A_2) \in \mathbb{R}^3$. Lastly, we define the projection $\mathcal{P}_0: H^3_{-\nu} (\mathbb{R}) \to \mathcal{E}_0$ by
\[\mathcal{P}_0 (q) = q_0 + q'(0) \xi + \frac{1}{2} q''(0) \xi^2, \]
noting that this is well-defined by Sobolev embedding and using that the solutions are actually in $H^k$ for $k$ as in the statement of Proposition \ref{p:cmfd}. 

The calculation of a reduced center flow is done in two steps: first, invariance is used to derive a Taylor expansion for the nonlocal center manifold $\Psi$ from Proposition \ref{p:cmfd}. Second, the flow $\Phi_\eta$ on the center manifold, which is defined by the action of translations  $\xi \mapsto \xi + \eta$ on $H^1_{-\eta}$, is projected onto the kernel and differentiated with respect to $\eta$ at $\eta = 0$, yielding a finite-dimensional reduced vector field. We note that the bulk of the computation in this process is in the first step, since the second step will consist entirely of differentiating polynomials. In fact, the entire computation involves only polynomials, highlighting the algebraic simplicity of the method.

In writing the center manifold as a graph $\Psi$ over $\mathcal{E}_0$, we seek a Taylor expansion of the form
\begin{equation}\label{e:tayexp} \Psi (A_0, A_1, A_2, \tbeta, \tc)= \sum_{\underset{|l|+|r|>1}{l, r}} A_0^{l_0} A_1^{l_1} A_2^{l_2} \tbeta^{r_1} \tc^{r_2} \psi_{l,r}(\xi), \end{equation}
  where $l = (l_1, l_2, l_3), r = (r_1, r_2)$. Here, the second multi-index $r$ is present because we are using the parameter-dependent version of \cite[Theorem 1]{fs2018center}. We will use invariance to solve for the Taylor coefficients $\psi_{l,r}$. Note that we do this using the unconditioned equation \eqref{e:FDEq}, since it makes calculations more straightforward and yields identical results.  We substitute \eqref{e:tayexp} into the functional differential equation \eqref{e:FDEq}, noting that 
\begin{equation}
    \mathcal{T}_{\tbeta,\tc}(q_0+\Psi) + \mathcal{N}(q_0+\Psi) = - \tc q_0'(\xi) - 2m \tbeta \big( q_0(\xi+1) - q_0(\xi) \big) +\mathcal{N}(q_0,0) + \mathcal{T}_{0,0}(\Psi) + \mathcal{O}(3),
\end{equation}
to obtain at quadratic order that
\begin{equation}\label{e:coeff1}
\begin{split}
     \sum_{|l| + |r| = 2} A_0^{l_0} A_1^{l_1} A_2^{l_2} \tbeta^{r_1} \tc^{r_2} \mathcal{T}_{0,0} \big(\psi_{l_0,l_1,l_2,r_1,r_2}(\xi) \big) &= \tc q_0'(\xi) + 2m \tbeta \big( q_0(\xi+1) + q_0(\xi) \big) -\mathcal{N}(q_0,0) \\ 
     &= 4 A_0A_1 + 4 \xi A_1^2 + 8 \xi A_0 A_2 + (12 \xi^2 + 4) A_1A_2 + (8 \xi^3 + 8 \xi + 4) A_2^2 \\ & \qquad - 2m A_1 \tbeta - A_1 \tc + (-2m-4m\xi)A_2 \tbeta  -2 \xi A_2 \tc.
     \end{split}
\end{equation}
We then take as an ansatz for $\psi_{l,r}$ the polynomials $ \psi_{l,r} = \sum_{i \ge 3} \alpha_i \xi^i$, suppressing the dependence on $l,r$ in $\alpha_i$. The ansatz is inspired by the fact that the kernel $\mathcal{E}_0$ consists of polynomials and more generally that the space of polynomials is invariant under convolution.
We calculate 
\begin{align}\label{e:coeff2} \mathcal{T}_{0,0}(\alpha_3 \xi^3 + \alpha_4 \xi^4 + \alpha_5 \xi^5 + \alpha_6 \xi^6) = (4m) \alpha_3 &+ (16m \xi - 24m)\alpha_4 + (40 m \xi^2 - 120m \xi +4m) \alpha_5 \nonumber
\\ & + (80\xi^3-360m \xi^2+24m \xi -120m) \alpha_6, \end{align}
and compare coefficients between \eqref{e:coeff2} and \eqref{e:coeff1} for each quadratic power of $(A_0, A_1, A_2, \tbeta, \wt{c})$. 
 After doing so, one finds that the nonzero Taylor coefficients at quadratic order are
\begin{equation}\label{e:psicoeff} \begin{split} &\psi_{1,1,0,0,0} = - \frac{1}{m} \xi^3, \hspace{0.5cm} \psi_{0,2,0,0,0} = - \frac{1}{4m} \xi^4 - \frac{3}{2m} \xi^3, \hspace{0.5cm} \psi_{1,0,1,0,0} = - \frac{1}{2m} \xi^4 - \frac{3}{m} \xi^3, \\
 &\psi_{0,1,1,0,0} = - \frac{3}{10m} \xi^5 - \frac{9}{4m} \xi^4 - \frac{71}{5m} \xi^3, \hspace{0.5cm} \psi_{0,0,2,0,0} = - \frac{1}{10m} \xi^6 - \frac{9}{10m}\xi^5 - \frac{71}{10m}\xi^4 - \frac{457}{10m} \xi^3, \\
& \psi_{0,1,0,1,0} =  \frac{1}{2} \xi^3, \hspace{0.5cm} \psi_{0,1,0,0,1} =  \frac{1}{4m} \xi^3, \hspace{0.5cm}
\psi_{0,0,1,1,0} =  \frac{1}{4} \xi^4 + 2 \xi^3, \hspace{0.5cm} \psi_{0,0,1,0,1} =  \frac{1}{8m} \xi^4 + \frac{3}{4m} \xi^3.\end{split}\end{equation}
 We can now compute the reduced vector field on the center manifold. We start by noting that the flow on the center manifold is defined by the action of translations, that is, at finite order,
 \begin{equation}
\begin{split} \Phi_\eta((q_0 + \Psi)(\xi)) &= (q_0 + \Psi)(\xi + \eta)\\  &= q_0(\xi+ \eta) + \sum A_0^{l_0} A_1^{l_1} A_2^{l_2} \tbeta^{r_1} \tc^{r_2} \psi_{l,r}(\xi + \eta), \end{split} \end{equation}
with $q_0, \psi_{l,r}$ as in \eqref{e:coeff1}, \eqref{e:coeff2}. This flow becomes a reduced flow after projection by $\mathcal{P}_0$ onto the kernel, and lastly becomes a reduced vector field after differentiation at $\eta = 0$: 
\begin{equation}\label{e:redform}
   q_0' = \frac{d}{d\eta} \mathcal{P}_0(\Phi_\eta(q_0 + \Psi))\Big|_{\eta = 0}.
\end{equation}
In order to express \eqref{e:redform} in terms of $(A_0, A_1, A_2)$, we note that $(q_0 + \Psi)(\xi+\eta)$ is a sum of polynomials in $\xi+\eta$, and calculate
\begin{equation*}
    \frac{d}{d\eta}\mathcal{P}_0((\xi + \eta))\big|_{\eta=0} = (1,0,0), \frac{d}{d\eta}\mathcal{P}_0((\xi + \eta)^2)\big|_{\eta=0} = (0,2,0), \frac{d}{d\eta}\mathcal{P}_0((\xi + \eta)^3)\big|_{\eta=0} = (0,0,3), 
\end{equation*}
\[\frac{d}{d\eta}\mathcal{P}_0((\xi + \eta)^n)\big|_{\eta=0} = (0,0,0), n \neq 1,2,3.\]


Then, using the expressions for $\Psi$ and $q_0$ in \eqref{e:psicoeff} and \eqref{e:q_0} and gathering terms, we find the reduced vector field to be given by
\begin{equation}\label{e:cmred}
\begin{split}
    \frac{d A_0}{d \eta} &= A_1 + \mathcal{O}(3),\\
    \frac{d A_1}{d \eta} &= 2A_2 + \mathcal{O}(3),\\
    \frac{d A_2}{d \eta} &=  - \frac{3}{m} \left( A_0A_1 + \frac{3}{2} A_1^2 + 3 A_0A_2 - \frac{m}{2} A_1 \tbeta - \frac{1}{4} A_1 \tc -2m A_2 \tbeta - \frac{3}{4} A_2 \tc  + \frac{71}{5} A_1A_2 + \frac{457}{10} A_2^2\right) +\mathcal{O}(3).
    \end{split}
\end{equation}

\subsection{Existence of solutions to reduced equations}\label{ss:exred}
It turns out that this three-dimensional ODE possesses a conserved quantity at leading order. Within level sets of this conserved quantity, one finds at leading order a homoclinic solution. In order to prove persistence of the homoclinic, one needs to prevent drift along level sets under perturbations of arbitrarily high order. To this aim, we establish the existence of a quantity that is \emph{exactly} conserved, which we compute to leading order following the ideas in \cite{bakker2018hamiltonian}. Within level sets of this exact conserved quantity, we then use a somewhat standard Melnikov-type argument to find homoclinics for the perturbed equation. 


In order to derive a conserved quantity, we integrate \eqref{e:FDE} from $0$ to $L$, $L \in \R$, and simplify, to find  \begin{align*}
     0 = &- c Q(L) + \int_{L-1}^{L} Q(y-1) Q(y+1) dy  -  \int_{L}^{L+1}Q(y-1) Q(y+1) dy  + \beta \int_{L}^{L+1}  Q^2(y) dy\\
     &- \left( -c Q(0) + \int_{-1}^{0} Q(y-1) Q(y+1) dy -  \int_{0}^{1} Q(y-1)Q(y+1)dy  + \beta \int_{0}^{1} Q^2(y)  dy \right),
 \end{align*}
 for any $L \in \R$. Hence $\Phi[Q]$, defined by
\begin{equation}\label{e:phi}
    \Phi[Q] = - c Q(0) +  \int_{-1}^{0} Q(y-1) Q(y+1) dy  -  \int_{0}^{1}Q(y-1) Q(y+1) dy  + \beta \int_{0}^{1}  Q^2(y) dy,
\end{equation}
is translation invariant on solutions, i.e. $\Phi[Q] = \Phi[S_{\tau} Q]$, where $S_{\tau}$ is a translation operator. We therefore define  $\varphi \equiv \Phi[Q]$, as the equivalent of a true conserved quantity or first integral.

Then, letting $Q(\xi) =  m +  A_0 + A_1 \xi + A_2 \xi^2 + \Psi(A_0, A_1, A_2, \wt{\beta},\wt{c})$ be an element of the center manifold, we insert $Q$ into \eqref{e:phi} to obtain  
\begin{equation} \label{e:varphi}
\begin{split}
    \varphi[Q] = (-2m^2 - m^2 \tbeta - m \tc) &+ \frac{4m}{3} A_2 + 2A_0^2 - (2m \tbeta + \tc) A_0 - \left( 4m \tbeta+ \frac{3}{2} \tc \right)A_1 + 6 A_0 A_1  + \frac{284}{15}A_0A_2 \\ 
    &+ \frac{142}{15} A_1^2 - \left(\frac{187 m}{15} \tbeta + \frac{22}{5} \tc \right) A_2 + \frac{457}{5}A_1A_2 + \frac{50201}{175}A_2^2 + \mathcal{O}(3).
\end{split}
\end{equation} 
We see that there exists a locally invertible change of coordinates where $A_2$ maps to $\phi$. We can solve \eqref{e:varphi} for $A_2$ using the implicit function theorem, since $\frac{4m}{3} \neq 0$, to find
\begin{equation}
    \begin{split}
 A_2 = \frac{3}{4m} \phi \ +  &\frac{3}{4m}(2 m \tbeta + \tc ) A_0 - \frac{3}{2m}A_0^2 - \frac{9}{2m} A_0 A_1 + \frac{3}{4m}\left(4 m \beta + \frac{3}{2} \tc \right) A_1 - \frac{71}{10m}A_1^2 \\& + \frac{99}{40m^2} \tc \phi +  \frac{561}{80} \tbeta \phi - \frac{213}{20m^2} A_0 \phi - \frac{4113}{80 m^2} A_1 \phi - \frac{1355427}{11200m^3} \phi^2 + \mathcal{O}(3) =: g(A_0, A_1, \phi) + \mathcal{O}(3),\end{split} \end{equation}
where $\phi$ is redefined to equal $\varphi + 2m^2 + m^2 \tbeta + m \tc$, since the latter is also conserved. We also define $\wt{A}_1 = \frac{d A_0}{d\eta} = A_1 + \mathcal{O}(3)$, and note that by the implicit function theorem we can  find $A_1$ in terms of $\wt{A}_1$ with an expansion. Changing coordinates to $(A_0, \wt{A}_1, \phi)$, we have the system 
\begin{equation}
    \begin{split}
        \frac{d A_0}{d \eta} &= \wt{A}_1, \\
        \frac{d \wt{A}_1}{d \eta} &= 2 g(A_0, \wt{A}_1,\phi) + \mathcal{O}(3),\\
        \frac{d \phi}{d \eta} &=  0.
    \end{split}
\end{equation}
We now rescale. Let $\tc = \tbeta m c_0 $, and rescale the spatial variable by choosing $\eta = (3\tbeta)^{\frac{1}{2}}(1 + \frac{c_0}{2})^\half \zeta$, the amplitudes by $A_0 = \tbeta m(1+\frac{c_0}{2})a_0$, $\wt{A}_1 = \tbeta^{\frac{3}{2}} \sqrt{3} m(1+\frac{c_0}{2})^\frac{3}{2} a_1$, and $\phi = \tbeta^2 3m^2(1+\frac{c_0}{2})^2 \wt{\phi}$.  Substituting, we obtain the final reduced system
\begin{equation}\label{eq:syst}
    \begin{split}
        \frac{d a_0}{d \zeta} &= a_1, \\
        \frac{d a_1}{d \zeta} &= \frac{1}{2}\wt{\phi} + a_0 - a_0^2 -\tbeta^\half\sqrt{3} \left[3(1+\frac{c_0}{2})^\half a_0a_1 -(1+\frac{c_0}{2})^{-\half}(2+\frac{3c_0}{4})a_1 \right] + \mathcal{O}(\tbeta),\\
        \frac{d \wt{\phi}}{d \zeta} &=  0.
    \end{split}
\end{equation}
We are looking for homoclinics to $a_0=0$ since we are interested in homoclinics with background mass $m$. We therefore choose $\wt{\phi} \equiv 0$ as a solution to the third equation, reducing to 
\begin{equation}\label{e:2dred}
    0 = a_0'' -a_0 + a_0^2 -\tbeta^\half\sqrt{3} \left[3(1+\frac{c_0}{2})^\half a_0a_0'-(1+\frac{c_0}{2})^{-\half}(2+\frac{3c_0}{4})a_0' \right] + \mathcal{O}(\tbeta) =:F(a_0,\tbeta),
\end{equation}
where $' = \partial_\zeta$. Note that for $\tbeta = 0$, this equation admits the explicit homoclinic solution $a_*(\zeta) = \frac{3}{2}\sech^2(\frac{\zeta}{2})$. \\
We now prove existence of homoclinic solutions to \eqref{e:2dred} using a standard Melnikov analysis. 

\begin{Lemma}
There exists $\tbeta_* > 0$ such that for $0 < \tbeta < \tbeta_*$, \eqref{e:2dred} admits a homoclinic solution to $a_0=0$ that is uniformly $\mathcal{O}(\tbeta^{1/2})$-close to $a_*(\zeta)$.
\end{Lemma}

\begin{Proof} Let $a_*(\zeta)$ be as above. Linearizing \eqref{e:2dred} at $a_0=a_*, \tbeta = 0$ yields the linear operator
$$ \mathcal{L}_* u= \partial_{\zeta\zeta}u + 2a_*u - u. $$
We note that $\ker \mathcal{L}_* = $ span$\{a_*'\}$. Define the projection $\mathcal{P}$ onto $\ker \mathcal{L}^*$ by $\mathcal{P} u = \langle  a_*',u \rangle$, where $\langle \cdot, \cdot \rangle$ denotes the $L^2$-inner product over $\mathbb{R}$. Note also that $\mathcal{L}_*$ is a self-adjoint operator. Then, writing $a_0 = a_* + v + \alpha a_*'$, where $v \in (\ker\mathcal{L}_*)^\perp$, we wish to solve $F(a_* + v + \alpha a_*',\tbeta)=0$. Without loss of generality, we fix a translate of any potential solution by choosing $\alpha = 0$, arriving at
\begin{equation}
    \begin{split}
        0 = \mathcal{P}F(a_* + v,\tbeta) \\
        0 = (1 - \mathcal{P})F(a_* + v, \tbeta).
    \end{split}
\end{equation}
By the implicit function theorem, since the linearization $(1 - \mathcal{P})DF(a_*,\tbeta) = (1 - \mathcal{P})\mathcal{L}_*$ is invertible as an operator from $\ker\mathcal{L}^\perp$ to $\ran\mathcal{L}$,  there exists a smooth function $v = \psi(\tbeta)$ defined on a neighborhood of $0$ such that $(1 - \mathcal{P})F(a_* + \psi(\tbeta), \tbeta)=0$. Inserting $\psi$ into the first equation, we get the reduced equation 
\begin{equation}
\begin{split}
   0 &=  \mathcal{P}F(a_* + \psi(\tbeta),\tbeta) \\
    &= \big\langle a_*',\mathcal{L}_*\psi(v) \big\rangle + \tbeta^\half\sqrt{3}\left[3(1+\frac{c_0}{2})^\half\big\langle a_*', a_*a_*'\big\rangle -  (1+\frac{c_0}{2})^{-\half}(2+\frac{3c_0}{4})\big\langle a_*', a_*'\big\rangle \right] + \mathcal{O}(\tbeta).
   \end{split}
\end{equation}
Noting that $\big\langle a_*',\mathcal{L}_*\psi(v) \big\rangle = \big\langle \mathcal{L}_*a_*',\psi(v) \big\rangle = 0$, and dividing by $\tbeta^\half,$ we get 
\begin{equation}\label{e:doublyred}
  0  = 3(1+\frac{c_0}{2})\big\langle a_*', a_*a_*'\big\rangle -  (2+\frac{3c_0}{4})\big\langle a_*', a_*'\big\rangle  + \mathcal{O}(\tbeta^\half).
\end{equation}
At $\tbeta = 0$, we can explicitly find $c_0 = \frac{4\langle a_*',a_*'\rangle}{3\langle a_*'', a_*''\rangle} - 2 = -\frac{16}{15}$. Noting that \eqref{e:doublyred} is smooth in $\delta = \tbeta^\half$, there exists $c_0 = c_0(\delta)$ for $|\delta| < \delta_*$, by the implicit function theorem, since $\frac{\partial}{\partial c_0}\mathcal{P}F 
= \frac{3}{4} \langle a_*'',a_*''\rangle = \frac{18}{7} \neq 0$, such that \eqref{e:doublyred} is satisfied.

Then for $0 < \tbeta < \tbeta_* := \delta_*^2$ there exists a solution to \eqref{e:2dred}, with the scaled correction $c_0$ to the speed given by  
\[c_0(\beta) = -\frac{16}{15} + \mathcal{O}(\tbeta^\half).\] \end{Proof}

Putting this together, we see that Theorem 1 is proven, and we have that a homoclinic solution to the system \eqref{e:cmred} exists for $\beta \lesssim 2$, with the speed parameter given by
\begin{equation*}c(\beta) = m(4-c_0\tbeta) + \mathcal{O}(\tbeta^\frac{3}{2}) 
= \frac{4m}{15}\left(7+4\beta\right) + \mathcal{O}((2-\beta)^\frac{3}{2}) .\end{equation*}


\section{Party drift: numerical continuation and stability}\label{s:num}

We explore existence and stability of  traveling parties numerically for $0 < \beta < 2$. The numerical results connect the two asymptotics regimes $\beta \gtrsim 0$ (\S{} \ref{s:small}) and $\beta\lesssim 2$ (\S{}\ref{s:large}). In particular, the results confirm the asymptotics and existence results with good quantitative agreement, and provide a more global picture of drift in the parameter $\beta$. In addition to these numerical continuation studies, we present a glimpse into the intricate question of stability and selection: instabilities in the constant party tails grow and lead to the formation of new parties. 

\subsection{Connecting the regimes -- secant continuation}

We compute traveling profiles throughout the entire range $0 < \beta < 2$ using a Newton method to find fixed points of the functional equation \eqref{e:fdefull},
\begin{equation}\label{e:fdefull2}
        0 =  -c q'(\xi) + m(2q(\xi -1) + 2(\beta + 1)q(\xi + 1) - 2(\beta + 1)q(\xi) - q(\xi -2) - q(\xi+2) + \mathcal{N}(q,\beta)), \qquad \xi\in\R.
\end{equation}
 We truncate the real line $\xi\in\R$ to $\xi\in[-L/2,L/2]$ and impose (artificial) periodic boundary conditions. We then discretize \eqref{e:fdefull2} with $N$ points and a grid spacing $h=L/N=1/\ell$, $\ell\in\N$, so that shifted values can be evaluated on the grid, using a 4th-order finite difference approximation for the derivative.

The linearization at a given profile $q_*$ possesses a two-dimensional generalized kernel, at least, generated by $q_*'$ from translations and by $q_*$ from mass scaling. We therefore add constraints: $\int \xi q_*=0$  eliminates translations, and $\int q_*=M$ fixes the mass. We compensate for the lack of a Lagrange multiplier associated with the mass constraint through the introduction of a dummy mass loss term  $\mu q(\xi)$. In summary, we solve the system
\begin{align}\label{e:disc}
             -c q'(\xi) + m\big[2q(\xi -1) + 2(\beta + 1)\big(q(\xi + 1) - q(\xi)\big) - & q(\xi -2) - q(\xi+2) + \mathcal{N}(q,\beta)\big] + \mu q = 0, \ \  \xi \in [-\frac{L}{2},\frac{L}{2}], \\
\int_{-L/2}^{L/2} \xi q(\xi) d\xi &=0 ,\\
 \int_{-L/2}^{L/2} (q(\xi) - m)d\xi&=1,\label{e:discmass}\\
 q(-L/2)&=q(L/2),\label{e:disc2}
 \end{align}
 after discretization in $\xi$ for the $N+2$ variables $(q,c,\mu)$ using a Newton method. We then add the parameter $\beta$ as a variable and a standard secant condition to the system \eqref{e:disc}---\eqref{e:disc2} for numerical continuation in $\beta$. We consistently find $\mu=0$ to machine precision as expected. 

We implement secant continuation starting at $\beta = 0.3$, computing profiles and speeds as $\beta$ varies. In the small-$\beta$ regime, we hold the total mass constant, and, as $\beta$ approaches 2, we hold the mass of the uniform background state constant, appropriately changing \eqref{e:discmass}. The initial interval is of width $L = 10$ in $x$ with $N = 1030$ total grid points. The grid is refined adaptively as $\beta \to 0$ by doubling the number of grid points when the second derivative of $q$ exceeds $0.03$ in the sup norm, since the profiles develop corners at the peak and sides. As $\beta \nearrow 2$,  the profiles get increasingly wider. We double the width $L$ of the interval and the number of grid points whenever the number of opinions at distance $0.2 L$ from the center exceeds $0.01$ of the peak value. We checked relative discretization and truncation errors by reducing $h$ and increasing $L$, and found errors typically of the order $10^{-7}$, always bounded by $10^{-4}$. 

We see below in Figure \ref{fig:partyshape} three examples of profiles as $\beta$ varies. We also plot in Figure \ref{fig:wake} the background mass as a function of $\beta$ when fixing $\int(q-m_\infty)=1$, that is, fixing the net mass in the party relative to the constant background distribution of opinions. Note that as $\beta \nearrow 2$ this mass approaches infinity. As $\beta \to 0$ the background mass decays exponentially in $\frac{1}{\beta}$ and therefore is difficult to compute precisely due to limitations in tolerances for the Newton method. 


\begin{figure}[H]
    \centering
    \begin{subfigure}{0.49\textwidth}
        \centering
        \includegraphics[width=\textwidth]{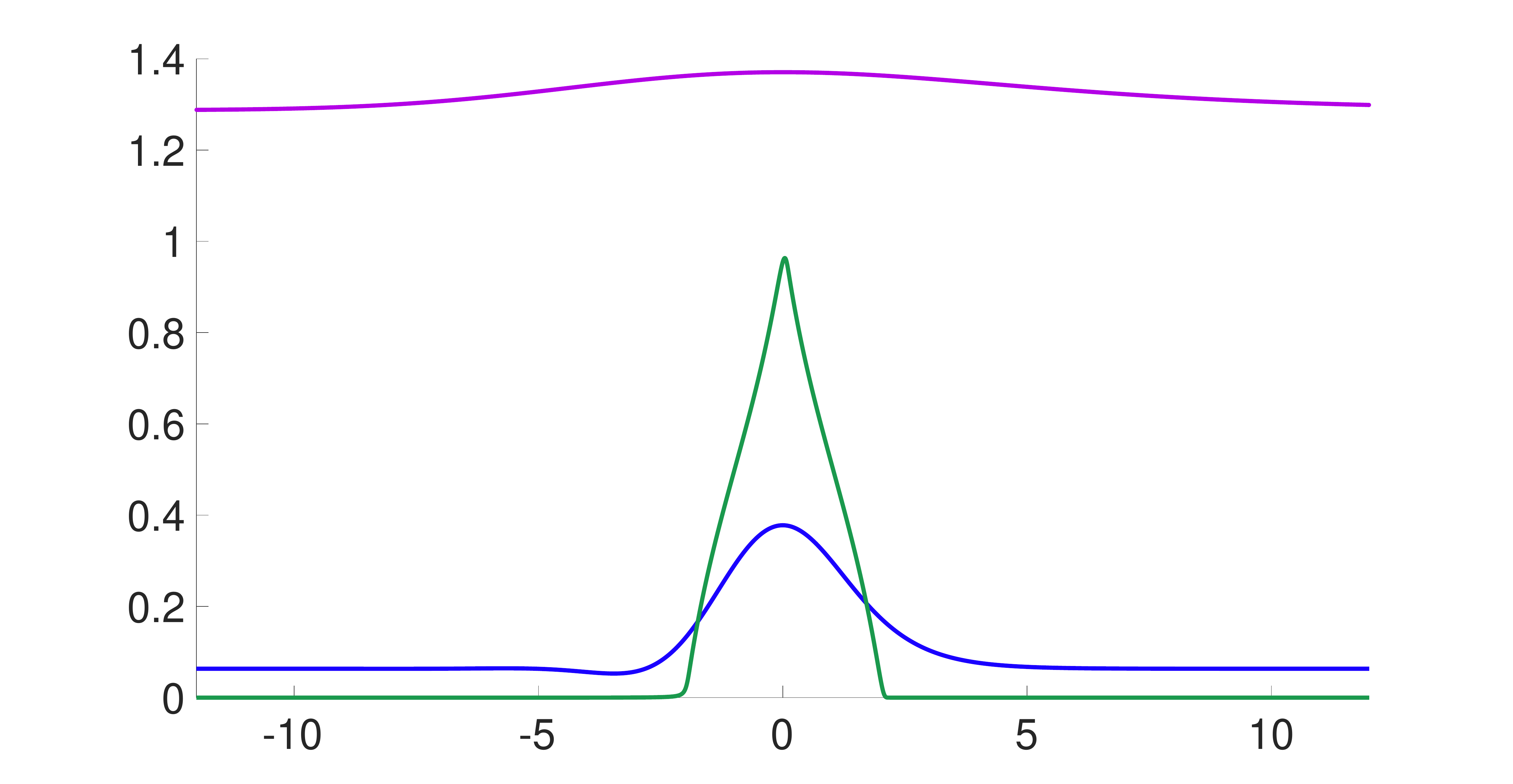}
        \caption{} 
        \label{fig:partyshape} 
    \end{subfigure}
        \begin{subfigure}{0.49\textwidth}
         \centering
         \includegraphics[width=\textwidth]{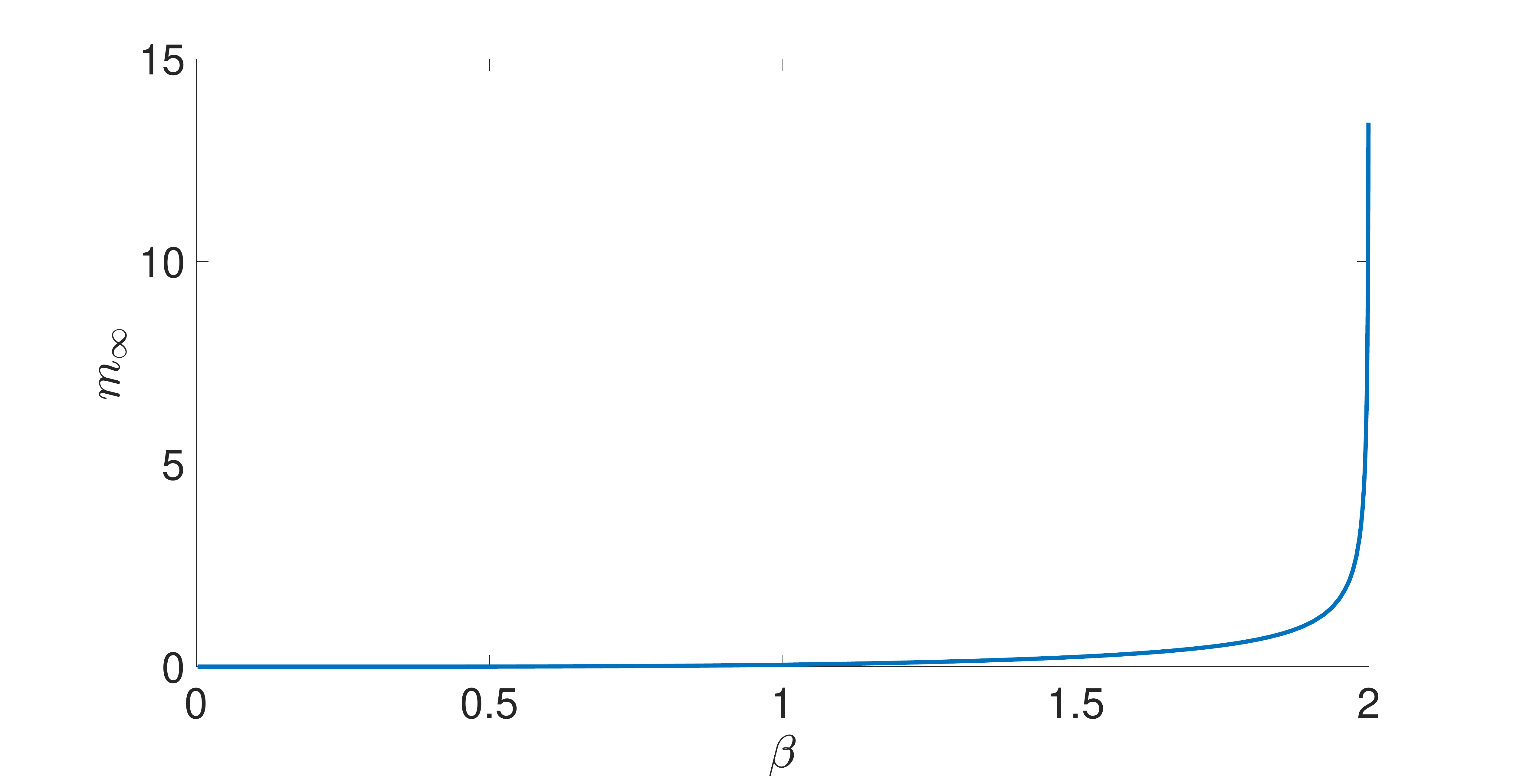}
         \caption{}
         \label{fig:wake}
     \end{subfigure}
    \caption{(a) Party shape for $\beta= .0006$ (green), $\beta= 1.08$ (blue), and $\beta = 1.92$ (purple), with net mass $m_{\textrm{party}} = 1$. (b) Plot of $\beta$ against the value of the uniform background mass $m_\infty$ needed to support a party of size 1.}
    \label{fig:cont}
\end{figure}

From the secant continuation one also obtains the speed $c$ of the profiles as $\beta$ varies. In Figure \ref{fig:coherent-party} below, we see the relationship between speed and $\beta$ normalized over total mass, party mass, and background mass, respectively. When $\beta$ is small, we fix the mass of the party at 1 and plot the speed, comparing it to its theoretical value 
from the drift speed calculations, including the numerically computed $\beta^{3/2}$ correction.  As $\beta$ approaches 2 we fix the uniform background mass and compare the speed with its theoretical value calculated in \S{}\ref{ss:exred}. In all cases we find excellent quantitative agreement with the leading-order predictions. 
\begin{figure}[H]
    \centering
    \begin{subfigure}{0.32\textwidth}
        \centering
         \includegraphics[width=\textwidth]{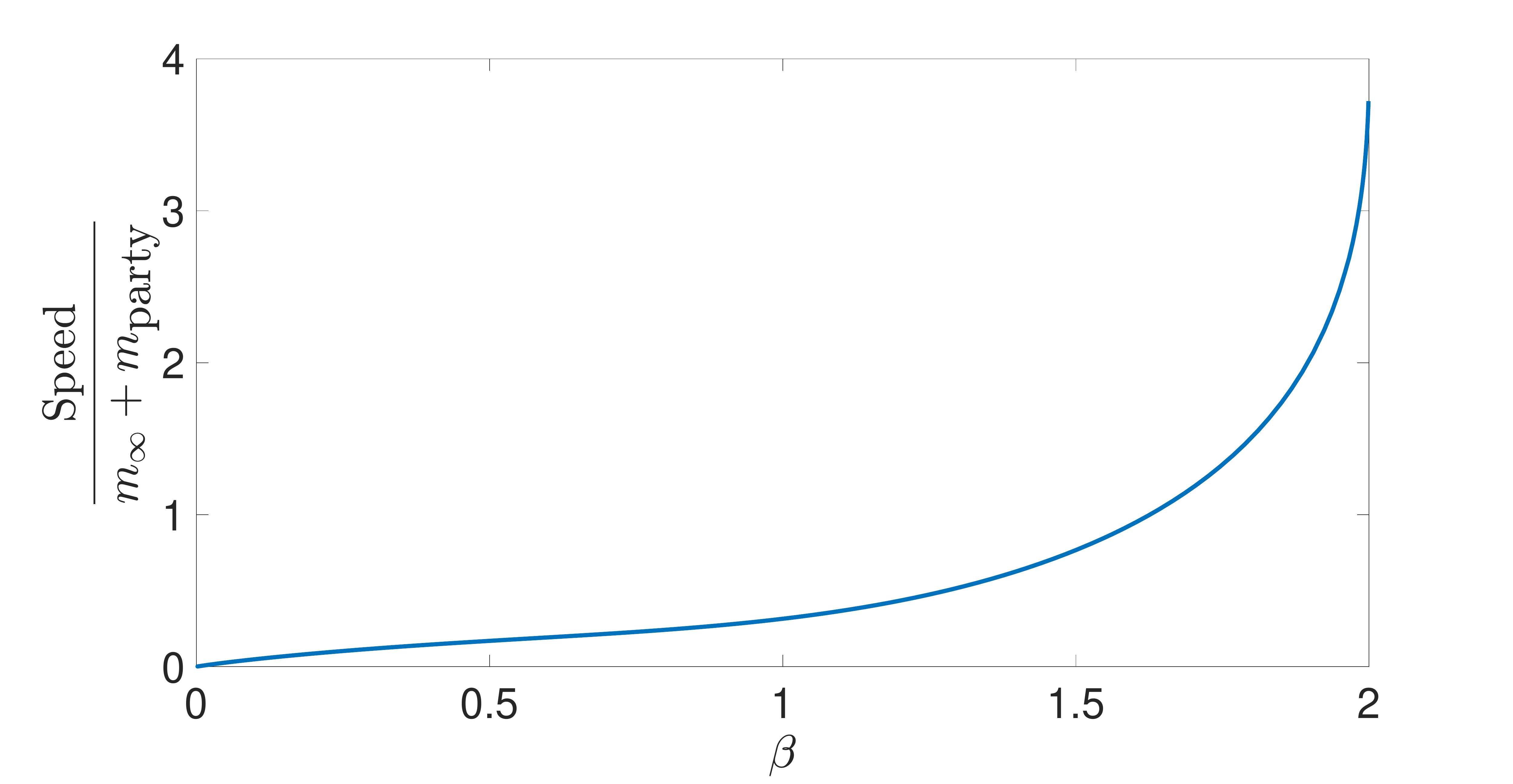}
         \caption{}
         \label{spreadingspeed} 
    \end{subfigure}
    \begin{subfigure}{0.32\textwidth}
         \centering
         \includegraphics[width=\textwidth]{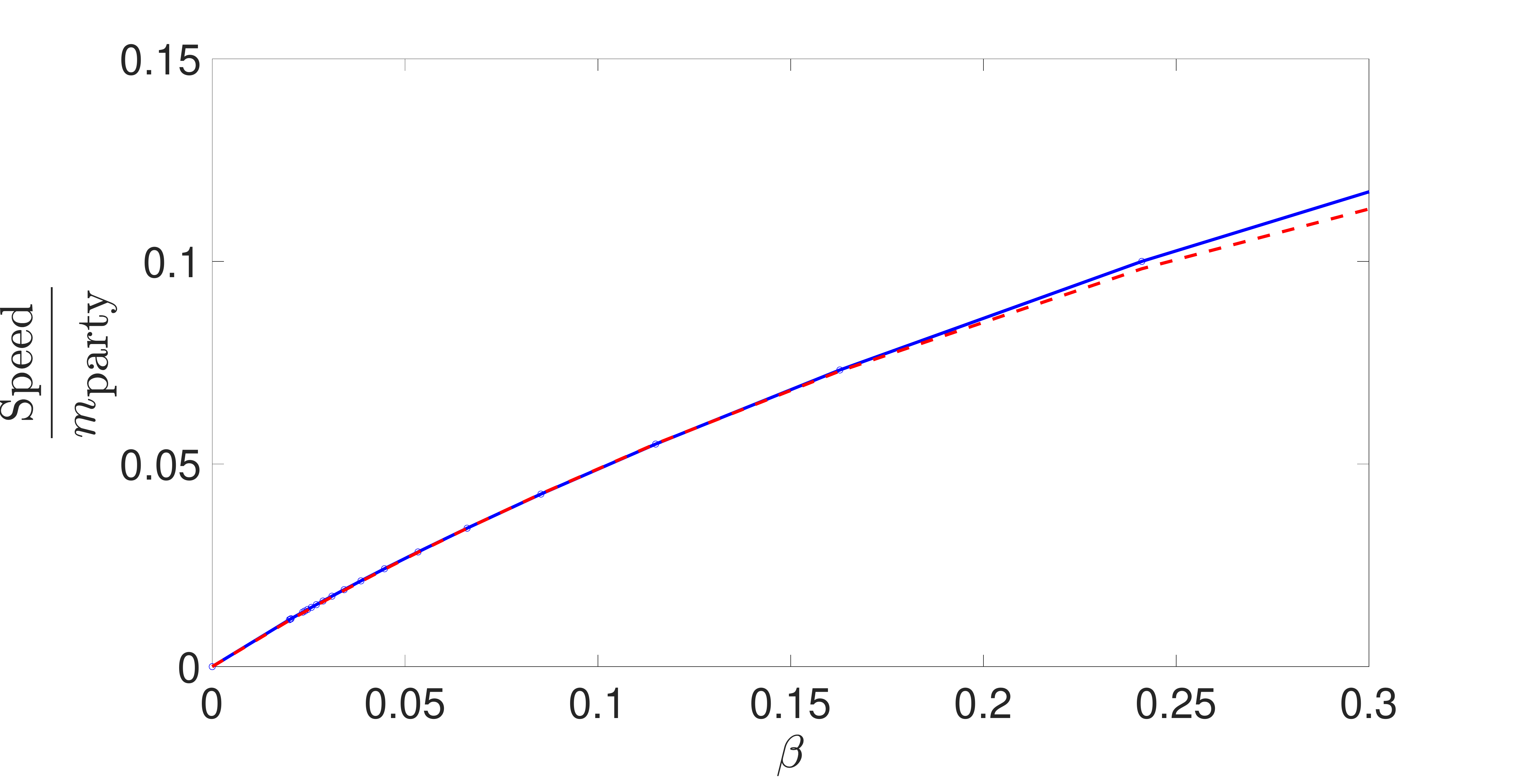}
         \caption{}
    \end{subfigure}
    \begin{subfigure}{0.32\textwidth}
         \centering
         \includegraphics[width=\textwidth]{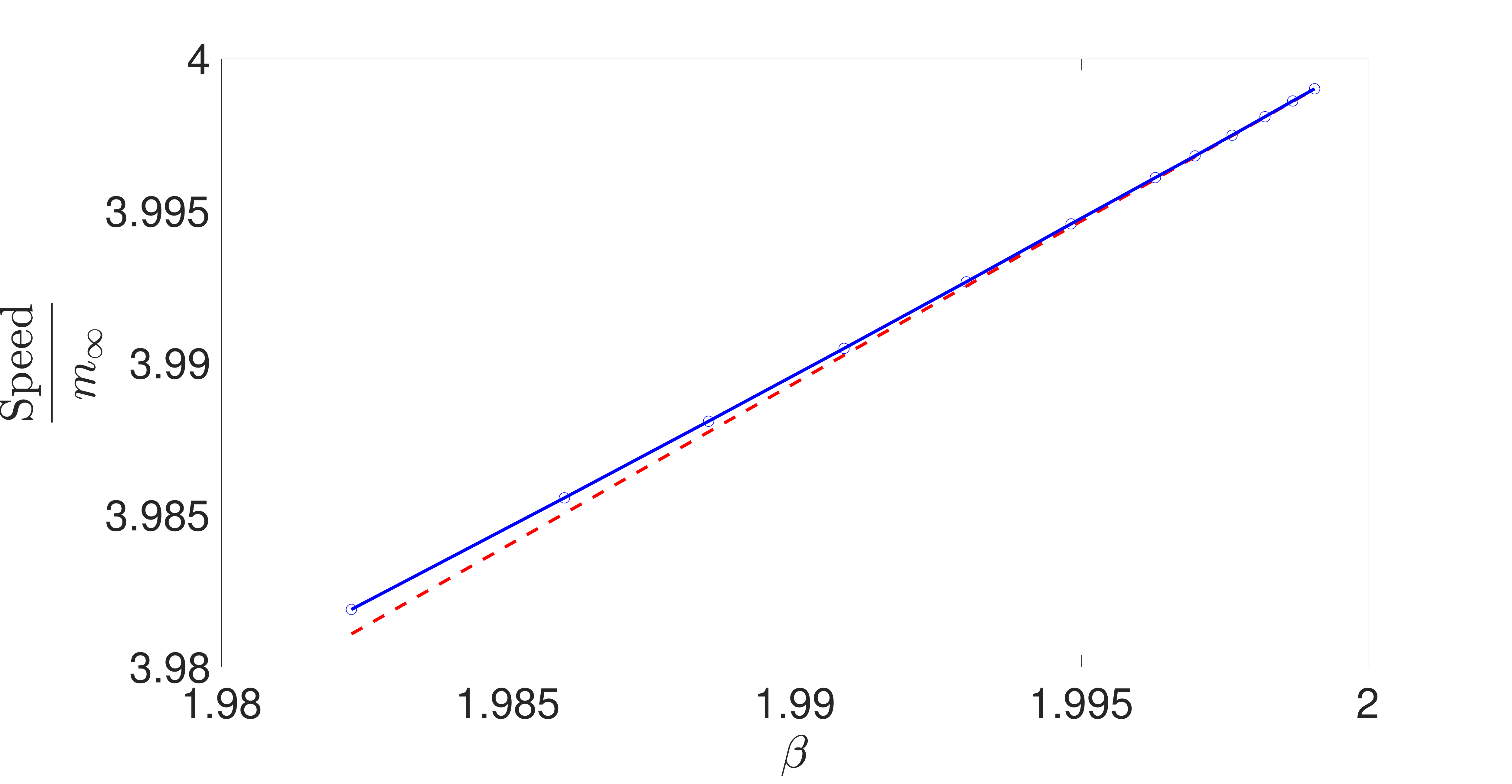}
         \caption{}
    \end{subfigure}
    \caption{Numerically computed speed for varying $\beta$ holding the sum of the background mass and party mass constant ((a)), holding the party mass constant ((b)),  
    and holding the background mass constant ((c)), with theoretical value shown in red in (b) and (c). 
    . }
    \label{fig:coherent-party}
\end{figure}





\subsection{Stability: direct simulations, convective stability, and pushed fronts}
Since  the uniform background state of a drifting party is unstable for $0 < \beta < 2$, drifting parties will typically not be stable in a strict sense. In fact, the instability of the background state is reflected in unstable continuous spectrum of the linearization. One then would wish to determine if the instability in the background affects the traveling party; that is, whether perturbations grow locally in a vicinity of the party, or whether perturbations are advected away from the center of the party, decaying locally uniformly while growing in norm. This distinction is commonly referred to as the difference between an absolute instability, where perturbations grow locally, and a convective instability, where perturbations decay locally uniformly in the frame comoving with the party. 

To a good approximation, this question is answered by comparing the spreading speed of perturbations of the uniform state $P_n=m_\infty,n\in\Z$ with the speed of the traveling party. Spreading speeds, in most scenarios, are determined as pulled or pushed speeds. In the case of pulled speeds, the spreading is determined by the linear equation, whereas for pushed speeds, the nonlinearity accelerates the propagation. We therefore start with the computation of the linear, pulled spreading speed; see \cite{VANSAARLOOS200329,MR3228472,bennaimscheel2016} for background and spreading in the non-biased case. In fact, spreading is mediated by pulled fronts, at the pulled, linear spreading speed, in the non-biased case as demonstrated in \cite{bennaimscheel2016}.

\begin{figure}[H]
    \centering
    \begin{subfigure}{0.49\textwidth}
         \centering
         \includegraphics[width=\textwidth]{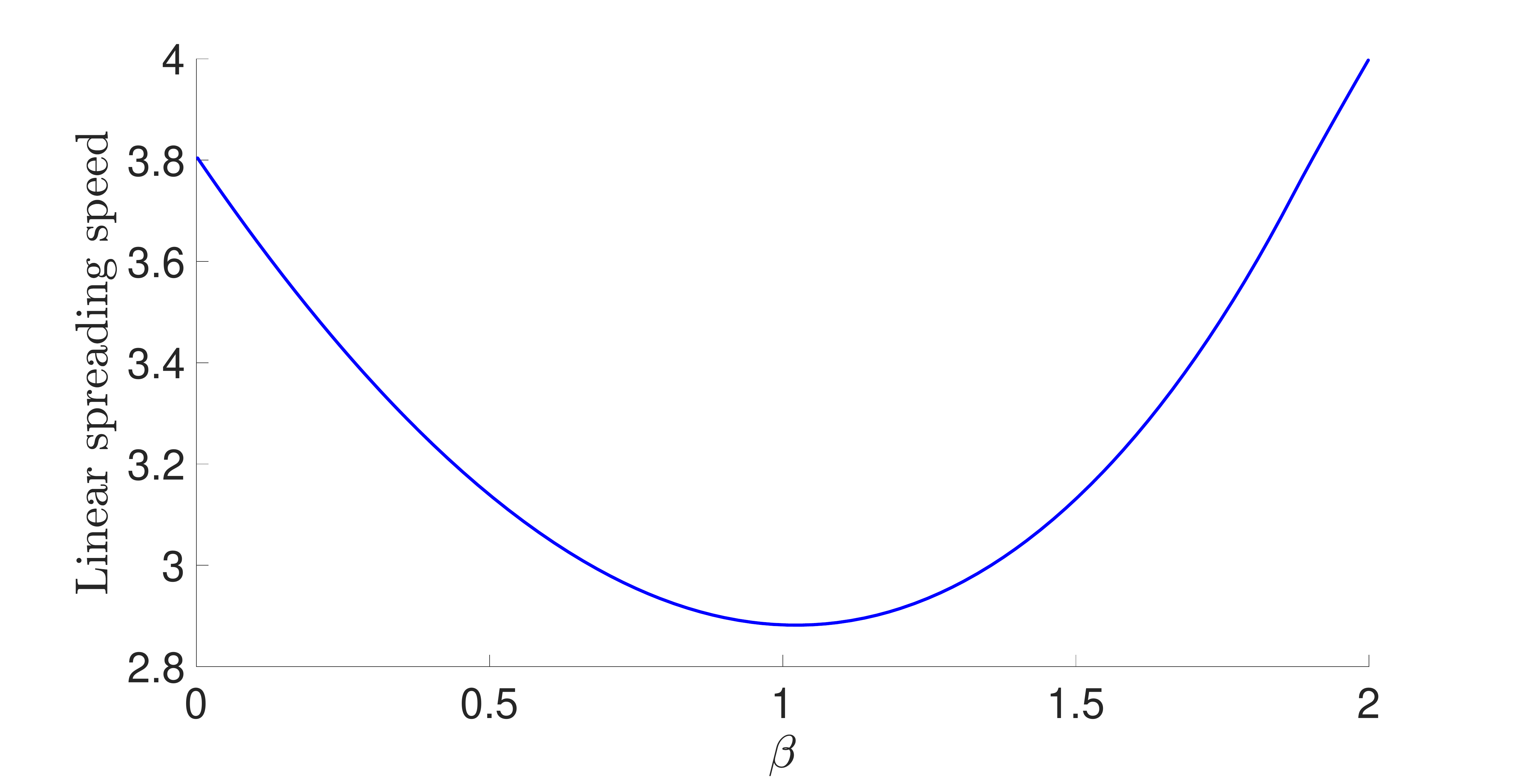}  
         \caption{}
         \label{fig:spreadingspeed}
     \end{subfigure}
     \begin{subfigure}{0.49\textwidth}
         \centering
\includegraphics[width=\textwidth]{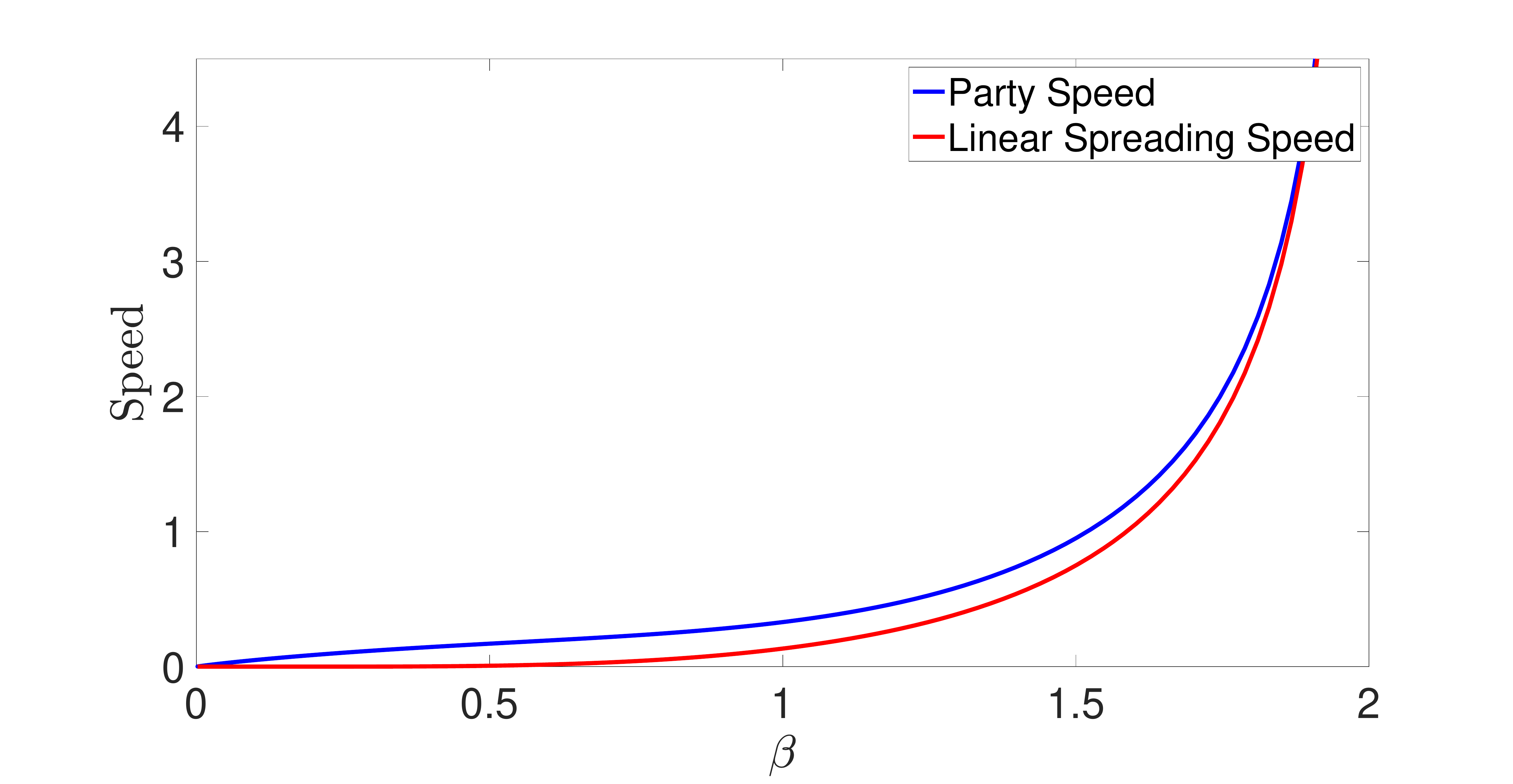} 
         \caption{}
         \label{fig:speedcomp}
     \end{subfigure}

    \caption{(a) Spreading speed of instabilities for $\beta \in [0,2]$ for $m=1$ (b) Party speed (blue) and linear spreading speed (red) for $\beta \in [0,2]$, with $m=m_{\infty}(\beta)$.  }
    \label{fig:spreading}
    
\end{figure}

Recall the dispersion relation for the linearized equation at the uniform state $P_j \equiv m$,
\begin{equation*}
    -i\omega = m(-2\cos(2\sigma) + (4+2\beta)\cos(\sigma) -(2+2\beta) + i\beta\sin(\sigma)). 
\end{equation*}
Then, through a saddle point analysis \cite{VANSAARLOOS200329,bennaimscheel2016}, the linear propagation velocity for any given $\beta$ is given by 
\begin{equation}
   v=\frac{d\omega}{d\sigma}=\frac{\Im[\omega]}{\Im[\sigma]},
   \label{eq:spreadspeed}
\end{equation}
after solving the second, complex equation for the complex variable $\sigma\in\C$. 
Figure \ref{fig:spreadingspeed} shows the numerical speeds $v$ with $\beta \in [0,2]$ and $m=1$, and Figure \ref{fig:speedcomp} compares this linear spreading speed for $m = m_\infty(\beta)$ with the party speed.
We find that the linear spreading speed of instabilities is less than the speed of the party for $\beta \in (0,2)$, which suggest that the party is linearly convectively unstable. 

\begin{figure}[H]
    \centering
        \includegraphics[width=.49\textwidth]{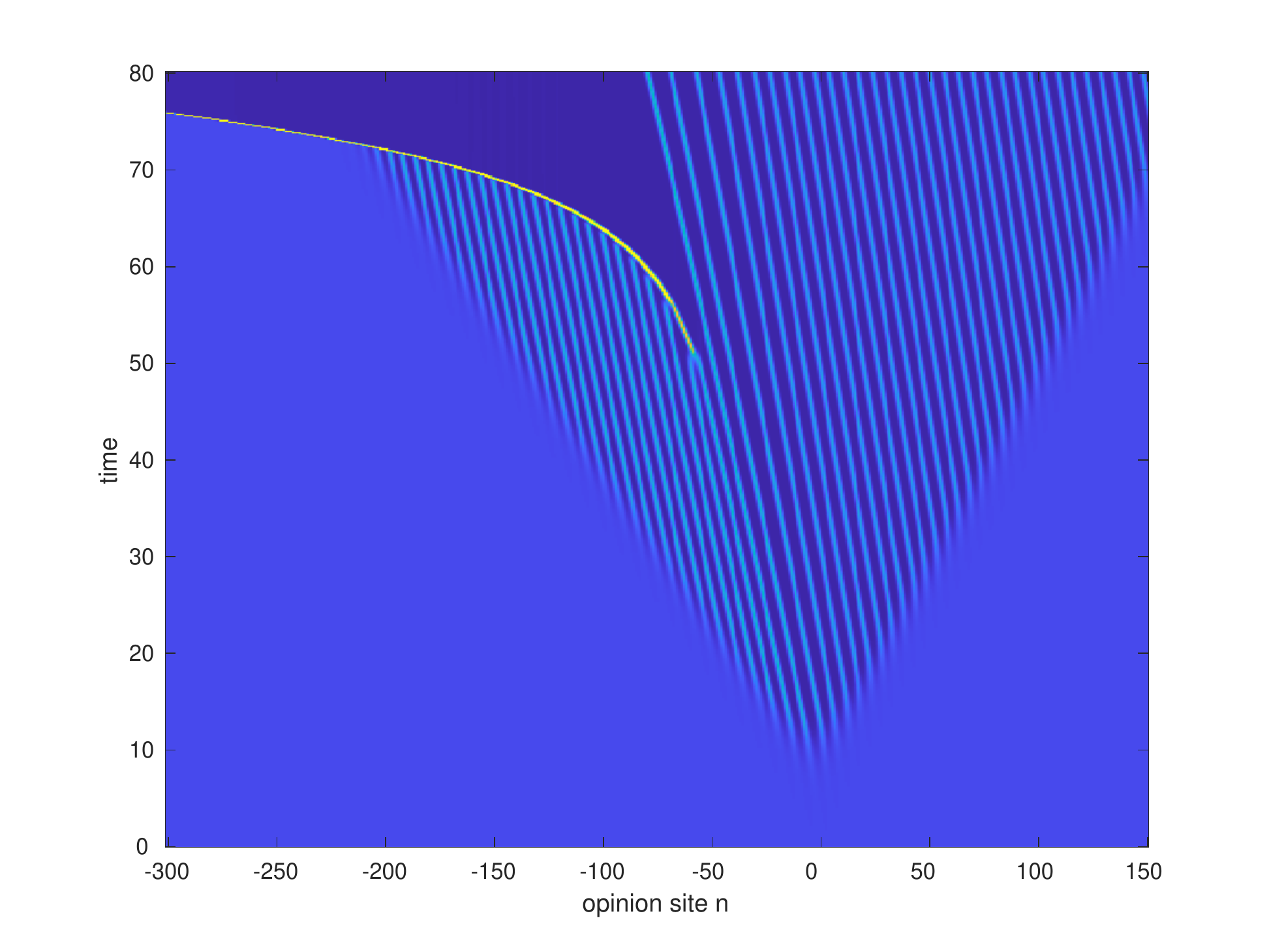}
         \ \ \includegraphics[width=.49\textwidth]{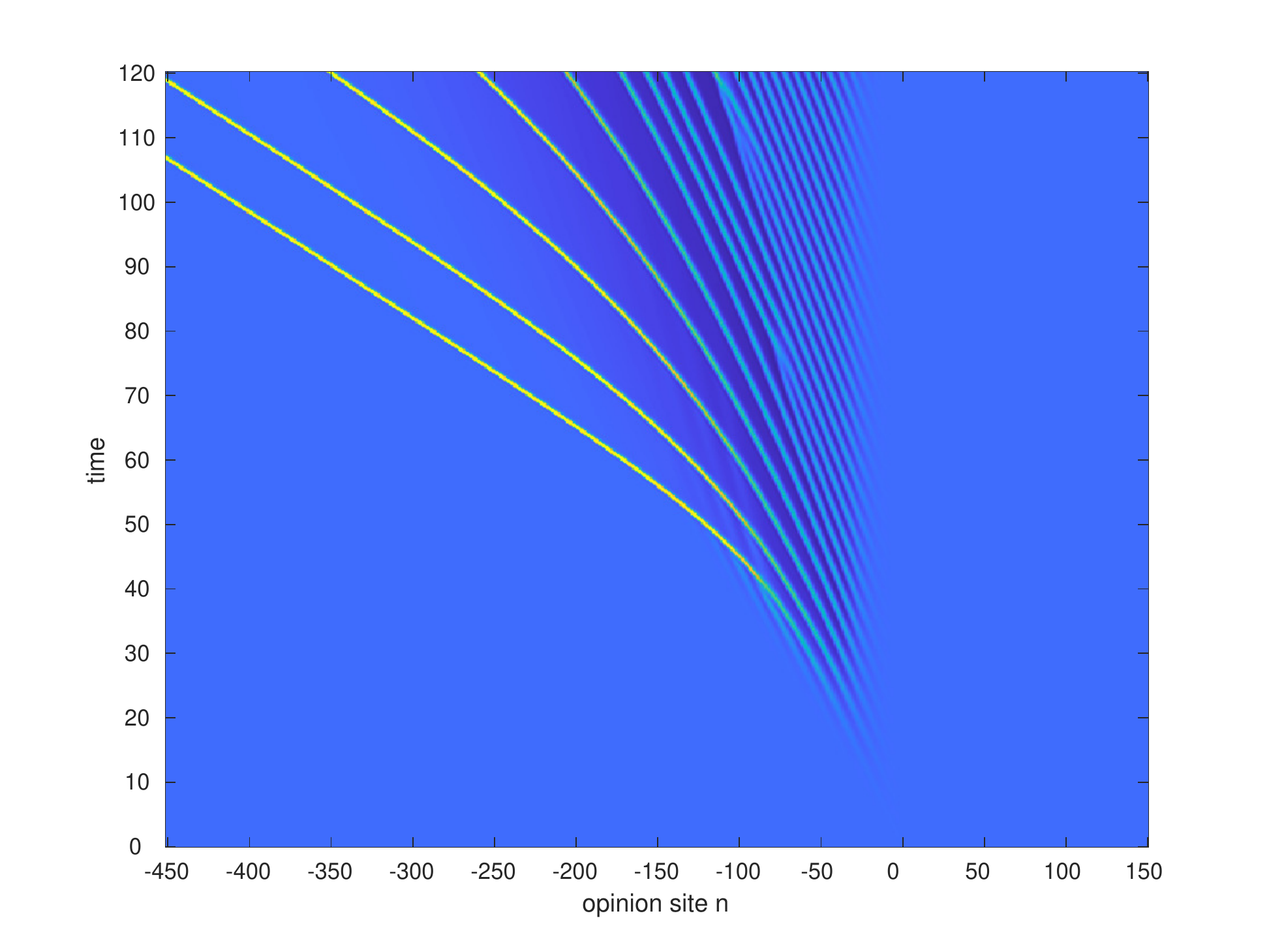}
    \caption{Space-time plots for $\beta= .35$ (left) and $\beta = 0.95$ (right), with the initial instability at $n=0$.}
    \label{fig:spacetime}
\end{figure}

We explore the possibility of a nonlinear instability, mediated by an analogue of a pushed front, through direct simulations. There, we do see evidence of faster than linear propagation. We simulate \eqref{e:bc_qb} in an adaptively moving frame in order to compute the spreading speed of instabilities, initializing at a small localized perturbation of the uniform state $P_j \equiv 1$. For all $\beta$, we see evidence of transient pulled fronts in the wake, which give way to sequences of parties traveling at a larger than linear speed. The initial transient of pulled front propagation is longer for smaller values of $\beta$. The effect is shown in  space-time plots (in a stationary frame) with initial conditions given as localized, small perturbations at $n=0$ of a uniform state $m = 1$. The transient pulled front is most clearly visible in the first plot, $\beta = 0.35,$ where we see a large single traveling party begin to form at $t = 50$ and overtake the pulled front at $t=75$. As $\beta$ increases, the transient pulled front has a narrow wake and is quickly reached by a large party that forms in its wake. In the right panel of Figure \ref{fig:spacetime}, the pulled front is visible from about $t = 10$ to $t = 50$. In this right panel, one also sees that after the initial nucleation of a large party that overtakes the pulled front, more large parties form subsequently spreading with similar speeds larger than the linear speed. The same effect occurs for smaller values of $\beta$, albeit on larger time scales. 
In conclusion, we suspect that the traveling parties are convectively stable, that is, observable, for long transients. The nonlinear process in the wake ultimately results in the creation of a train of traveling parties of comparable size, that remain well separated, traveling at similar speeds. 

A complete description of dynamics, either on a constant background, or on a zero background with resulting loss of mass, appears quite intricate. The space-time plots presented here suggest a similarity with phenomena observed in excitable media in simulations and experiments, coined trace-firing and backfiring \cite{doi:10.1063/1.2943307,doi:10.1063/1.2266993}, both of which are also poorly understood mathematically.





\section{Non-coherence for other bias terms}\label{s:nonincite}

Different bias mechanisms can generate dynamics very different from the dynamics in equation  \eqref{e:bc_qb} that we analyzed thus far.  We present here briefly phenomena caused by two other forms of bias, namely  bias in the compromise process and linear bias mimicking a transport term. For the former, propagation is blocked at one-site parties; for the latter, we find diffusive dissipation of the party, similar to the dynamics of the bounded confidence model with diffusion studied in \cite{Ben_Naim_2005}. We think of those two scenarios as evidence that the  self-incitement mechanism of \eqref{e:bc_qb} is  rather special in allowing  coherent movement of parties. 

\subsection{Bias in the compromise process}
A natural mechanism for introducing bias would be in the compromise process itself; that is, to consider an equation such as
\begin{equation}\label{e:biasincomp}
    \frac{d P_n}{d t} = (2-\beta)P_{n+1}P_{n-1} - (1-\beta)P_n P_{n+2} - P_n P_{n-2}),
\end{equation}
where agents interact in the same way as in \eqref{e:bc}, but the probability of changing the opinion to the compromise opinion is not equal for the two interacting agents.  Another formulation is to include bias in nearest-neighbor interactions and consider the equation
\begin{equation}\label{e:nonlinbias}
    \frac{dP_n}{dt} = 2P_{n+1}P_{n-1} - P_n(P_{n+2} + P_{n-2}) + \beta(P_{n+1}P_n-P_nP_{n-1}),
\end{equation}
where nearest-neighbor interactions lead to agents moving to the left with probability $\beta$. 

In both cases \eqref{e:biasincomp} and \eqref{e:nonlinbias}, the bias does not induce persistent drift of existing parties. In fact, one- and two-site parties are also equilibria of \eqref{e:biasincomp}, and one-site parties are equilibria of \eqref{e:nonlinbias}. As a result, parties do not drift at all in \eqref{e:biasincomp}. In \eqref{e:nonlinbias}, one can mimic the analysis in \S\ref{s:small} and find a drift speed near two-site parties, which vanishes at one-site parties. The result agrees well with direct numerical simulations which find two site parties evolving towards one-site parties, which are one-sided stable, similar to saddle-node equilibria. One can think of this blocking of motion at one-site parties as a pinning phenomenon, reflecting the discreteness of the opinion space, similar to the rather well understood pinning of front and pulse propagation in discrete or inhomogeneous media \cite{ankney2019}. 

We note however that the dynamics in \eqref{e:nonlinbias} resulting from the perturbation of a spatially constant equilibrium with the resulting formation of multiple parties is rather complex and does involve movement of parties, but at non-constant speeds, mediated by both the bias effect on single parties and the interaction between parties.


\subsection{Non-incitement bias}
All bias terms considered thus far are quadratic, modeling person-to-person interactions, and  preserve the quadratic scaling invariance of the dynamics. The arguably simplest possible bias terms would model spontaneous change of opinion in one direction, without need for interaction between agents, and thus be represented by a linear term of the form  $\beta(P_{n+1} - P_n)$, with variants $ \beta(P_{n+\ell} - P_n)/\ell$, $\ell=2,3,\ldots$. Up to scaling, this term can be viewed as a spatial discretization of the shift term $u_t = u_x$ in a continuous opinion space $x\in\R$. The resulting equation is
\begin{equation}\label{e:linbias}
    \frac{dP_n}{dt} = 2P_{n+1}P_{n-1} - P_n(P_{n+2} + P_{n-2}) + \beta(P_{n+\ell}-P_n)/\ell,
\end{equation} 
with several possible scalings depending on initial conditions. Numerically, we observe drift of single parties as expected, with initial speed $\beta$ at leading order. Expanding the one-sided differentce $P_{n+1}-P_n$ into derivatives, we find at second order an effective diffusion, which indeed is observable in the dynamics and leads to mass loss in the party similar to the observations in  \cite{Ben_Naim_2005}. Figure \ref{fig:tailmass} below shows the loss of mass from the party, defined as all opinions within 5 sites of the peak, into the tail and leading edge over time. Figure \ref{fig:frontbacktails} shows the self-similar shape of the leading edge when scaled horizontally by $t^{-\half}$,  mimicking the profiles in \cite{Ben_Naim_2005}. The tail, however, appears to grow exponentially in time at first, then becomes diffusive as the mass of the party no longer is significant; see also Figure \ref{fig:tailmass}. The effect of mass loss can also be understood in the small-$\beta$ perturbation analysis of \S\ref{s:small}, where at higher order, the drift term generates growth at sites with distance 3 or more to the party. In fact, for long-range coupling $\ell>2$, the effect appears at first order in $\beta$ and the mass loss is correspondingly more pronounced. 
\begin{figure}[H]
    \centering
    \begin{subfigure}{0.49\textwidth}
        \centering
        \includegraphics[width=\textwidth]{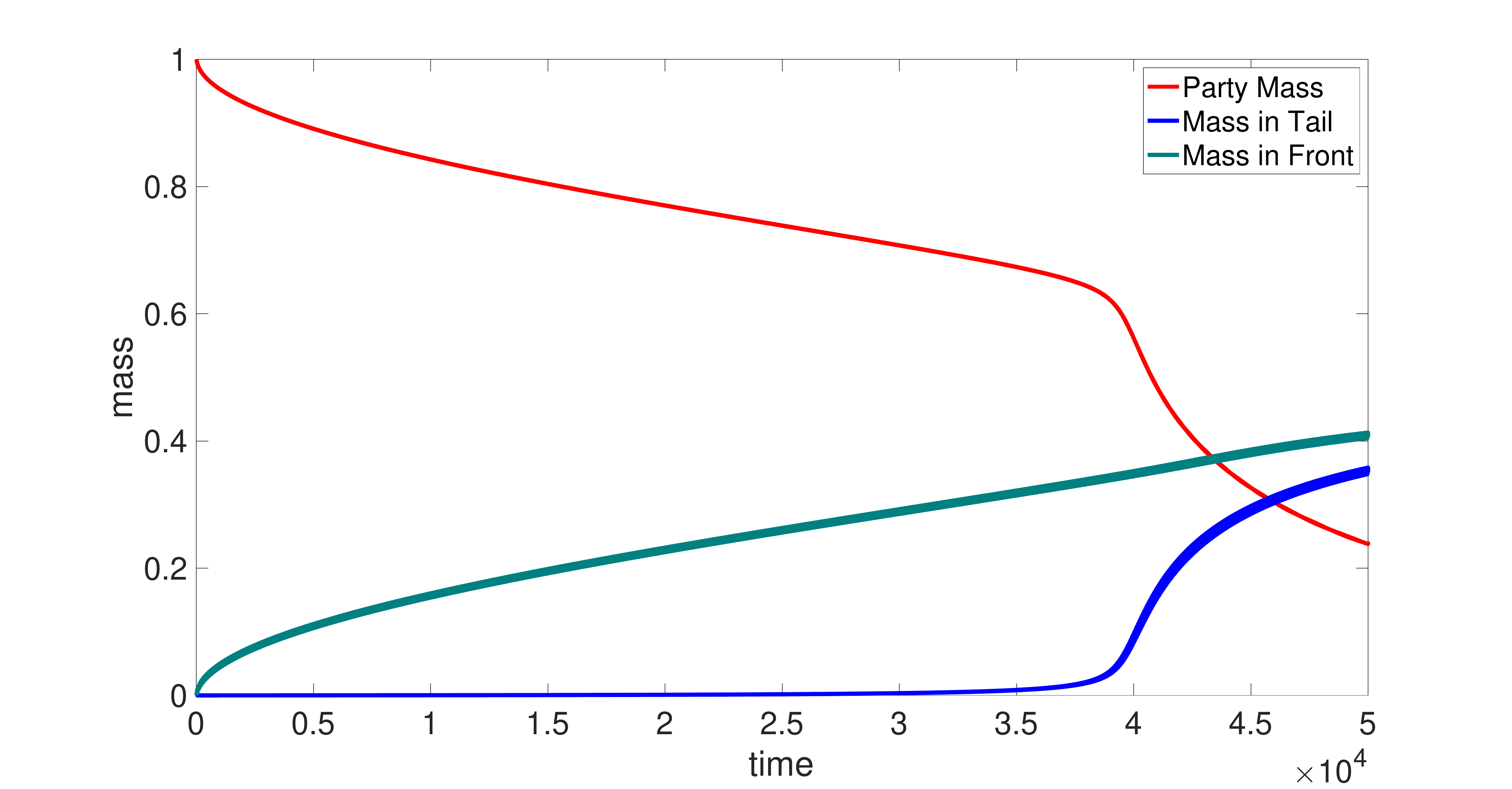}
        \caption{Proportion of mass in party, tail, and leading tail over time starting from a one-site cluster.} 
        \label{fig:tailmass} 
    \end{subfigure}
        \begin{subfigure}{0.49\textwidth}
         \centering
         \includegraphics[width=\textwidth]{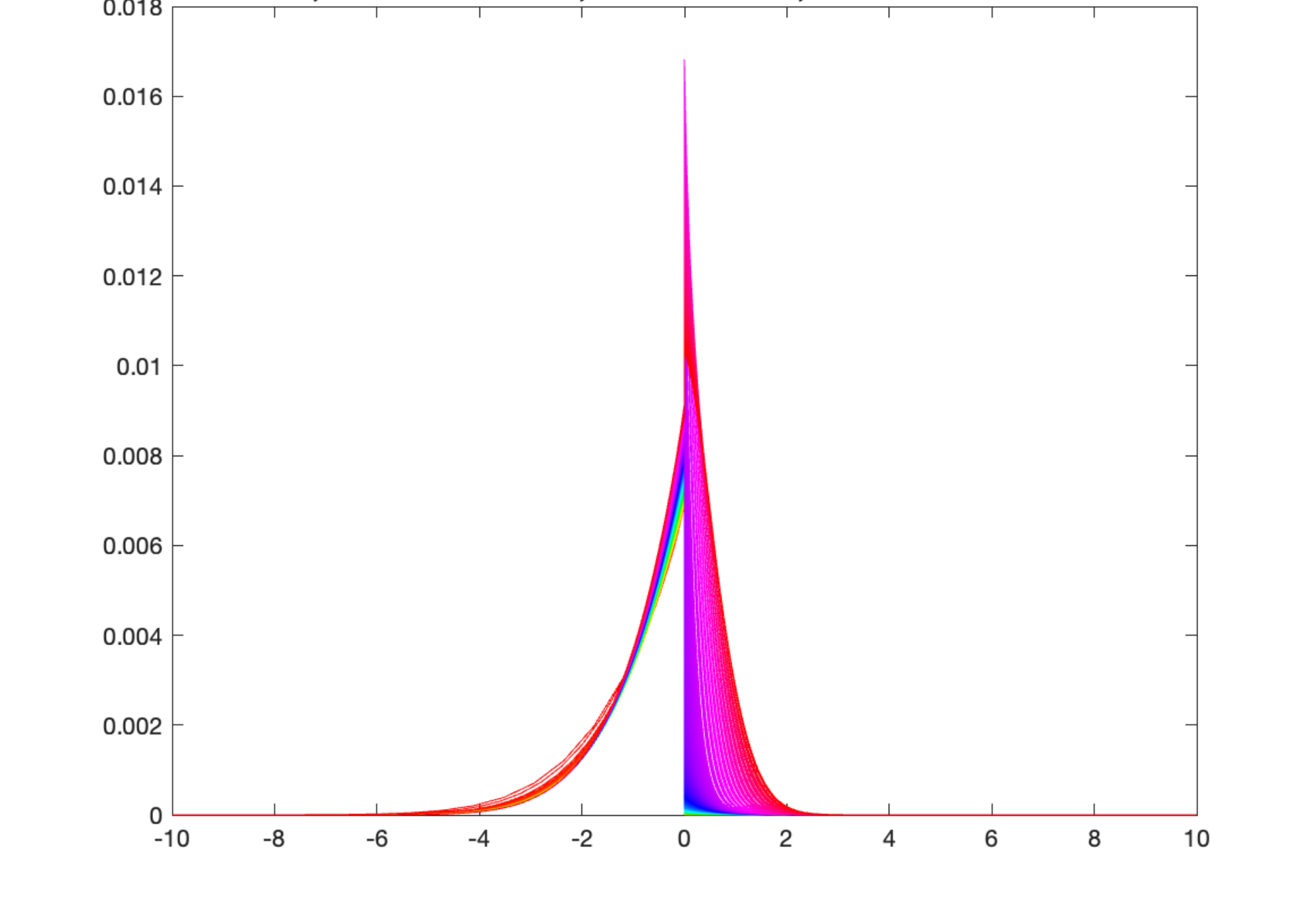}
         \caption{Profiles of the leading and trailing tails, omitting the party in between, scaled by plotting against $xt^\half$.}
         \label{fig:frontbacktails}
     \end{subfigure}
    \caption{Plot of proportion of mass in the party, tail, and leading tail over time, and profiles of the tails plotted against $t^\half x$}. 
    \label{fig:tails}
\end{figure}




\section{Discussion}

We presented results on the effect of bias in the bounded confidence model. Roughly speaking, the bounded confidence model supports localized clusters, including a one-parameter family roughly parameterized by the position of the cluster in space. One expects that introduction of bias leads to a drifting movement of clusters. We analyze this drifting motion in two limits, small and large bias. Notably, we prove rigorously that coherent movement of parties on a constant background distribution of opinions is possible for quadratic self-incitement bias in the large bias regime $\beta \lesssim 2$. Numerically, we find such coherent motion for all values of bias $0<\beta<2$. For small values of $\beta$, the constant background distribution is exponentially small in the parameter. 

Technically, we used a geometric singular perturbation analysis to derive drift speeds in the small bias regime and a nonlocal center manifold analysis to find coherent drifting parties for $\beta\lesssim 2$. The analysis in the latter case is possibly of independent interest, demonstrating the simplicity of calculations in the recently introduced framework of center manifolds without a phase space. Within the narrow focus on coherent drifting parties, a major open question is to establish rigorously existence, drift speeds, and size of background state in the regime $\beta\gtrsim 0$. 

Beyond coherent drifting parties, we touched on the evolution near unstable constant states. We find modulation equations familiar from fluid dynamics, such as the Korteweg-deVries or the Kuramoto-Sivashinsky equation and variations thereof. Observed dynamics in simulations appear to involve complex dynamics of formation and interaction of parties. 

Interesting questions arise when attempting to quantify mass loss in incoherent drifting parties. Geometric singular perturbation theory at a single party predicts mass accumulation in sites further away from the center of the party, possibly at high or beyond all orders in $\beta$ depending on the bias term. It would be interesting to quantify these effects and, for quadratic self-incitement bias, to contrast with the existence of coherent drifting parties on a constant background, relating for instance the size of the constant background to the rate of mass loss. 

We noted that stability questions in the context of bias are subtle, due to the instability of the constant background state and the complexity of the dynamics in the evolution of perturbations. In fact, it appears that even without bias, stability of parties against perturbations, say small in $\ell^\infty$ is not known.







\end{document}